\documentclass[11pt,a4paper]{article}
\usepackage{cite}
\usepackage{amsmath}
\usepackage{amssymb}

\usepackage{graphicx}
\usepackage{epsfig}

\newtheorem{corollary}{Corollary}
\newtheorem{lemma}{Lemma}
\newtheorem{proposition}{Proposition}
\newtheorem{definition}{Definition}

\newtheorem{theorem}{Theorem}
\newcommand{\onetom}{1,2,\cdots,m}

\usepackage[colorlinks]{hyperref}

\begin{document}

\title{Synchronization of discrete-time dynamical
networks with time-varying couplings}

\author{Wenlian Lu\thanks{School of
Mathematical Sciences, Fudan University, 200433, Shanghai, China
({\tt wenlian.lu@gmail.com}).} 
\and Fatihcan M. Atay\thanks{Max
Planck Institute for Mathematics in the Sciences, Inselstr.~22,
04103 Leipzig, Germany ({\tt atay@member.ams.org}).}
\and J\"urgen Jost\thanks{Max Planck Institute for Mathematics
in the Sciences, Inselstr.~22, 04013 Leipzig, Germany
({\tt jjost@mis.mpg.de}).}}

\date{Preprint. Final version in:\\
\textit{SIAM Journal on Mathematical Analysis}, 39(4):1231-1259, 2007\\
\href{http://dx.doi.org/10.1137/060657935}{http://dx.doi.org/10.1137/060657935}	}

\maketitle

\begin{abstract}
We study the local complete synchronization of discrete-time
dynamical networks with time-varying couplings. Our conditions for
the temporal variation of the  couplings are rather general and
include both variations in the network structure and in the
reaction dynamics; the reactions could, for example, be driven by
a random dynamical system. A basic tool is the concept of Hajnal
diameter which we extend to infinite Jacobian matrix sequences.
The Hajnal diameter can be used to verify synchronization and we
show that it is equivalent to other quantities which have been
extended to time-varying cases, such as the projection radius,
projection Lyapunov exponents, and transverse Lyapunov exponents.
Furthermore, these results are used to investigate the
synchronization problem in coupled map networks with time-varying
topologies and possibly directed and weighted edges. In this
case, the Hajnal diameter of the infinite coupling matrices can be
used to measure the synchronizability of the network process. As
we show, the network is capable of synchronizing some chaotic map
if and only if there exists an integer $T>0$ such that for any
time interval of length $T$, there exists a vertex which can
access other vertices by directed paths in that time interval.
\end{abstract}



\section{Introduction}
Synchronization of dynamical processes on networks is presently an
active research topic. It represents a mathematical framework that
on the one hand can elucidate -- desired or undesired --
synchronization phenomena in diverse applications. On the other
hand, the synchronization paradigm is formulated in such a manner
that powerful mathematical techniques from dynamical systems and
graph theory can be utilized. A standard version is
\begin{equation}
x^{i}(t+1)=f^{i}(x^{1}(t),x^{2}(t),\cdots,x^{m}(t)),\quad
i=\onetom, \label{simple}
\end{equation}
where $t\in{\mathbb Z}^{+}=\{0,1,2,\cdots\}$ denotes the discrete
time, $x^{i}(t)\in \mathbb{R}$ denotes the state variable of unit
(veretex) $i$, and for $i=\onetom$,
$f^{i}:\mathbb{R}^{m}\rightarrow \mathbb{R}$ is a $C^{1}$
function. This dynamical systems formulation contains two aspects.
One of them is the reaction dynamics at each node or vertex of the
network. The other one is the coupling structure, that is, whether
and how strongly, the dynamics at one node is directly influenced
by the states of the other nodes.

Equation (\ref{simple}) clearly is an abstraction and simplification
of synchronization problems found in applications. On the basis of
understanding the dynamics of (\ref{simple}), research should then
move on to more realistic scenarios. Therefore, in the present
work, we address the question of synchronization when the right
hand side of (\ref{simple}) is allowed to vary in time. Thus, not
only the dynamics itself is a temporal process, but also the
underlying structure changes in time, albeit in some
applications that may occur on a slower time scale.

The essence of the hypotheses on $f=[f^1,\cdots,f^m]$ needed for
synchronization results (to be stated in precise terms shortly) is
that synchronization is possible as an invariant state, that is,
when the dynamics starts on the diagonal $[x,\cdots,x]$, it will
stay there, and that this diagonal possesses
 a stable attracting  state. The question about synchronization then
 is whether this state is also attracting for dynamical states
$[x^1,\cdots,x^m]$ outside the diagonal, at least locally, that is
when the components $x^i$ are not necessarily equal, but close to
each other. This can be translated into a question about
transverse Lyapunov exponents, and one typically concludes that
the existence of a synchronized attractor in the sense of Milnor.
In our contribution, we can already strengthen this result by
concluding (under appropriate assumptions) the existence of a
synchronized attractor in the strong sense instead of only in the
weaker sense of Milnor. (We shall call this local complete
synchronization.) This comes about because we achieve a
reformulation of the synchronization problem in terms of Hajnal
diameters (a concept to be explained below).

Our work, however, goes beyond that. As already indicated, our
main contribution is that we can study the local complete
synchronization of general coupled networks with time-varying
coupling functions, in which each unit is dynamically evolving
according to
\begin{eqnarray}
x^{i}(t+1)=f^{i}_{t}(x^{1}(t),x^{2}(t),\cdots,x^{m}(t)),\quad
i=\onetom.\label{general}
\end{eqnarray}
This formulation, in fact, covers both aspects described above, the
reaction dynamics as well as the coupling structure. The main purpose
of the present paper then is to identify general conditions under which we can
prove synchronization of the dynamics (\ref{general}). Thus, we can
handle variations of the reaction dynamics as well as of the
underlying network topology. We shall mention below various
applications where this is of interest.

Before that, however, we state our technical hypothesis on the
right hand side of (\ref{general}): for each $t\in\mathbb Z^{+}$,
$f^{i}_{t}:\mathbb{R}^{m}\rightarrow \mathbb{R}$ is a $C^{1}$
function with the following hypothesis:

{\em$\bf H_{1}$. There exists a $C^{1}$ function
$f(s):\mathbb{R}\rightarrow \mathbb{R}$ such that
\begin{eqnarray*}
f^{i}_{t}(s,s,\cdots,s)=f(s)
\end{eqnarray*}
holds for all $s\in \mathbb{R}$, $t\in {\mathbb Z}^{+}$, and
$i=\onetom$. Moreover, for any compact set $K\subset\mathbb R^{m}$
$f^{i}_{t}$ and the Jacobian matrices $[\partial f^{i}_{t}/\partial
x^{j}]_{i,j=1}^{m}$ are all equicontinuous in $K$ with respect to
$t\in\mathbb Z^{+}$ and the latter are all nonsingular in $K$.}

This hypothesis ensures that the diagonal synchronization manifold
$${\mathcal S}=\bigg\{[x^{1},x^{2},\cdots,x^{m}]^{\top}\in\mathbb
R^{m}:~x^{i}=x^{j},~i,j=\onetom\bigg\}$$
 is an invariant manifold for
the evolution (\ref{general}). If
$x^{1}(t)=x^{2}(t)=\cdots=x^{m}(t)=s(t)$ denotes the synchronized
state, then
\begin{eqnarray}
s(t+1)=f(s(t)).\label{syn}
\end{eqnarray}
For the synchronized state (\ref{syn}), we assume the existence of
an attractor:

{\em $\bf H_{2}$. There exists a compact asymptotically stable
attractor $A$ for Eq.~(\ref{syn}). That is, (i)
$A\subset\mathbb{R}$ is a forward invariant set; (ii) for any
neighborhood $U$ of $A$ there exists a neighborhood $V$ of $A$
such that $f^{n}(V)\subset U$ for all $n\in\mathbb Z^{+}$; (iii)
for any sufficiently small neighborhood $U$ of $A$, $f^{n}(U)$
converges to $A$, in the sense that for any neighborhood $V$,
there exists $n_{0}$ such that $f^{n}(U)\subset V$ for $n\ge
n_{0}$; (iv) there exists $s^{*}\in A$ for which the
$\omega$-limit set is $A$. }

Let $A^m$ denote the Cartesian product $A\times\cdots\times A$
($m$ times).  Local complete synchronization (synchronization for
simplicity) is defined in the sense that the set $\mathcal{S}\cap
A^m=\{[x,\cdots ,x]: x \in A\}$ is an asymptotically stable
attractor in $\mathbb R^{m}$. That is, for the coupled dynamical
system (\ref{general}), differences between components converge to
zero if the initial states are picked sufficiently near
$\mathcal{S}\cap A^m$, i.e., if  the components are all close to
the attractor $A$ and if their differences are sufficiently small.
In order to show such a synchronization, one needs a third
hypothesis {\em $\bf H_{3}$} that in technical terms is about
Lyapunov exponents transverse to the diagonal. That is, while the
dynamics on the attractor may well be expanding (the attractor
might be chaotic), the transverse directions need to be suitably
contracting to ensure synchronization. The corresponding
hypothesis {\em $\bf H_{3}$} will be stated below (see
(\ref{thm1})) because it requires the introduction of crucial
technical concepts.

It is an important aspect of our work that we shall derive the
attractivity here in the classical sense, and not in the sense of
Milnor, i.e., not only some set of positive measure, but a full
neighborhood is attracted. For details about the difference
between Milnor attractors and asymptotically stable attractors;
see \cite{Ash,Mil1}. Usually, when studying synchronization, one
derives only the existence of a Milnor attractor; see \cite{Pik}.

The motivation for studying (\ref{general}) comes from the
well-known coupled map lattices (CML) \cite{Kaneko}, which can be written as
follows:
\begin{eqnarray}
x^{i}(t+1)=f(x^{i}(t))+\sum\limits_{j=1}^{m}L_{ij}f(x^{j}(t)),~i=\onetom,\label{cml-comp}
\end{eqnarray}
where $f:\mathbb{R}\rightarrow \mathbb{R}$ is a differentiable map
and $L=[L_{ij}]_{i,j=1}^{m}\in \mathbb{R}^{m\times m}$ is the
diffusion matrix, which is determined by the topological structure
of the network and satisfies $L_{ij}\ge 0$ for all $i\ne j$, and
$\sum_{j=1}^{m}L_{ij}=0$ for all $i=\onetom$. Letting $
x=[x^{1},x^{2}, \dots,x^{m}]^{\top}\in \mathbb{R}^{m}$, $
F(x)=[f(x^{1}),f(x^{2}),\dots,f(x^{m})]^{\top}\in \mathbb{R}^{m}$,
and $G=I_{m}+L$, where $I_{m}$ denotes the identity matrix of
dimension $m$, the CML (\ref{cml-comp}) can be written in the
matrix form
\begin{eqnarray}
x(t+1)=G F(x(t))\label{CMLs}
\end{eqnarray}
where $G=[G_{ij}]_{i,j=1}^{m}\in \mathbb{R}^{m\times m}$ denotes
the coupling and satisfies $G_{ij}\ge 0$ for $i\ne j$ and
$\sum_{j=1}^{m}G_{ij}=1$ for all $i=\onetom$. So, if $G_{ii}\ge 0$
holds for all $i=\onetom$, then $G$ is a stochastic matrix.

Recently, synchronization of CML has attracted increasing
attention \cite{Pik,Jost,Ding,Atay,Lu}. Linear stability analysis
of the synchronization manifold was proposed and transverse
Lyapunov exponents were used to analyze the influence of the
topological structure of networks. In \cite{Ash},
 conditions for generalized transverse stability were presented. If the
transverse (normal) Lyapunov exponents are negative, a chaotic
attractor on an invariant submanifold can be asymptotically stable
over the manifold.  Ref.~\cite{Wu,Lu1} have found out that chaos
synchronization in a network of nonlinear continuous-time or
discrete-time dynamical systems respectively is possible if and only
if the corresponding graph has a spanning tree. However,
synchronization analysis has so far been limited to autonomous
systems, where the interactions between the vertices (state
components) are static and do not vary through time.

In the social, natural, and engineering real-world, the topology
of the network often varies through time. In communication
networks, for example, one must consider dynamical networks of
moving agents. Since the agents are moving, some of the existing
connections can fail simply due to occurrence of an obstacle
between agents. Also, some new connections may be created when one
agent enters the effective region of other agents \cite{Sab}. On
top of that, randomness may also occur. In a communication
network, the information channel of two agents at each time may be
random \cite{Hat}. When an error occurs at some time, the
connections in the system will vary. In \cite{Sab,Hat,Vic},
synchronization of multi-agent networks was considered where the
state of each vertex is adapted according to the states of its
connected neighbors with switching connecting topologies. This
multi-agent dynamical network can be written in discrete-time form
as
\begin{eqnarray}
x^{i}(t+1)=\sum\limits_{j=1}^{m}G_{ij}(t)x^{j}(t),~i=\onetom,\label{multi}
\end{eqnarray}
where $x^{j}(t)\in\mathbb R$ is the state variable of vertex $j$
and $[G_{ij}(t)]_{i,j=1}^{m}$, $t\in\mathbb Z^{+}$, are stochastic
matrices. Ref.~\cite{Mor} considered a convexity-conserving coupling
function which is equivalent to the linear coupling function in
(\ref{multi}). It was found that the connectivity of the switching
graphs plays a key role in the synchronization of multi-agent
networks with switching topologies. Also, in the recent literature
\cite{Lv,Bel,Sti}, synchronization of continuous-time dynamical
networks with time-varying topologies was studied. The
time-varying couplings  investigated, however, are specific, with
either symmetry \cite{Lv},  node balance \cite{Bel}, or fixed time
average \cite{Sti}.

Therefore, it is natural to investigate the synchronization of CML
with general time-varying connections as:
\begin{eqnarray}
x(t+1)=G(t)F(x(t))\label{TVCMLs}
\end{eqnarray}
where $G(t)=[G_{ij}(t)]_{i,j=1}^{m}\in\mathbb R^{m\times m}$
denotes the coupling matrix at time $t$ and
$F(x)=[f(x_{1}),\cdots,f(x_{n})]^{\top}$ is a differentiable
function.
We shall address this problem in the context of the general
coupled system (\ref{general}).

Let
\begin{eqnarray*}
x(t)=\left[\begin{array}{c}x^{1}(t)\\x^{2}(t)\\\vdots\\x^{m}(t)\end{array}\right]
~\text{ and }~F_{t}(x(t))=\left[\begin{array}{c}
f^{1}_{t}(x^{1}(t),\cdots,x^{m}(t))\\f^{2}_{t}(x^{1}(t),\cdots,x^{m}(t))\\\vdots
\\f^{m}_{t}(x^{1}(t),\cdots,x^{m}(t))\end{array}\right].
\end{eqnarray*}
Eq. (\ref{general}) can be rewritten in a matrix form:
\begin{eqnarray}
x(t+1)=F_{t}(x(t)).\label{matrix}
\end{eqnarray}

The time-varying coupling can have a special form and may be
driven by some other dynamical system. Let $\mathcal
Y=\{\Omega,{\mathcal F},P,\theta^{(t)}\}$ denote a metric
dynamical system (MDS), where $\Omega$ is the metric state space,
${\mathcal F}$ is the $\sigma$-algebra, $P$ is the probability
measure, and $\theta^{(t)}$ is a semiflow satisfying
$\theta^{(t+s)}=\theta^{(t)}\circ\theta^{(s)}$ and
$\theta^{(0)}=\rm{id}$, where $\rm{id}$ denotes the identity map.
Then, the coupled system can  be regarded as a random dynamical
system (RDS) driven by $\mathcal Y$:
\begin{eqnarray}
x(t+1)=F(x(t),\theta^{(t)}\omega),~t\in\mathbb
Z^{+},~\omega\in\Omega.\label{RDS}
\end{eqnarray}
In fact, one can regard the dynamical system (\ref{RDS}) as a skew
product semiflow,
\begin{eqnarray*}
&&\Theta: \mathbb Z^{+}\times \Omega\times \mathbb
R^{m}\rightarrow\Omega\times\mathbb R^{m}\\
&&\Theta^{(t)}(\omega,x)=(\theta^{(t)}\omega,x(t)).
\end{eqnarray*}
Furthermore, the coupled system can have the form
\begin{eqnarray}
x(t+1)=F(x(t),u(t)),~t\in\mathbb Z^{+},\label{ex}
\end{eqnarray}
where $u$ belongs to some function class $\mathcal U$ and may be
interpreted as an external input or force. Then, defining
$[\theta^{(t)}u](\tau)=u(t+\tau)$ as a shift map, the system
(\ref{ex}) has the form of (\ref{RDS}). In this paper, we first
investigate the general time-varying case of the system
(\ref{matrix}) and also apply our results to systems of the form
(\ref{RDS}).

To study synchronization of the system (\ref{matrix}), we use its
variational equation by linearizing it. Consider the difference
$\delta x^{i}(t)=x^{i}(t)-f^{(t-t_{0})}(s_{0})$. This implies that
$\delta x^{i}(t)-\delta x^{j}(t)=x^{i}(t)-x^{j}(t)$ holds for all
$i,j=\onetom$. We have
\begin{eqnarray}
\delta x^{i}(t+t_{0})=\sum\limits_{j=1}^{m}\frac{\partial
f^{i}_{t+t_{0}-1}}{\partial x^{j}}(f^{(t-1)}(s_{0}))\delta
x^{j}(t+t_{0}-1),~i=\onetom.\label{var-comp}
\end{eqnarray}
where for simplicity we have used the notation
 $\frac{\partial f^{i}_{t+t_{0}-1}}{\partial
x^{j}}(f^{(t-1)}(s_{0}))$ to denote\\ $\frac{\partial
f^{i}_{t+t_{0}-1}}{\partial
x^{j}}(f^{(t-1)}(s_{0}),\cdots,f^{(t-1)}(s_{0}))$. Let
$$\delta x(t)=\left[\begin{array}{c}\delta
x^{1}(t)\\\vdots\\\delta x^{m}(t)\end{array}\right],\quad
D_{t}(s)=\bigg[\frac{\partial f^{i}_{t}}{\partial
x^{j}}(s)\bigg]_{i,j=1}^{m}.$$ The variational equation
(\ref{var-comp}) is written in matrix form,
\begin{eqnarray}
\delta x(t+t_{0})=D_{t+t_{0}-1}(f^{(t-1)}(s_{0}))\delta
x(t+t_{0}-1).\label{var-matrix}
\end{eqnarray}

For the Jacobian matrix, the following lemma is an immediate
consequence of the hypothesis $\bf{H_{1}}$.
\begin{lemma}
\label{lem1.1}
\begin{eqnarray*}
\sum\limits_{j=1}^{m}\frac{\partial f^{i}_{t}}{\partial
x_{j}}(s,s,\cdots,s)=f'(s),\quad i=\onetom~\text{ and
}~t\in{\mathbb Z}^{+}.
\end{eqnarray*}
\end{lemma}
Namely, all rows of the Jacobian matrix $[\partial
f^{i}_{t}/\partial x_{j}]_{i,j=1}^{m}$ evaluated on the
synchronization manifold ${\mathcal S}$ have the same sum, which
is equal to $f'(s)$.

As a special case, if the time variation is driven by some
dynamical system $\mathcal Y=\{\Omega,\mathcal
F,P,\theta^{(t)}\}$, then the variational system does not depend
on the initial time $t_{0}$, but only on $(s_{0},\omega)$. Thus, the
Jacobian matrix can be written in the form
$D(f^{(t)}(s_{0}),\theta^{(t)}\omega)=D_{t}(f^{(t)}(s))$, by which
the variational system can be written as:
\begin{eqnarray}
\delta x(t+1)=D(f^{(t)}(s_{0}),\theta^{(t)}\omega)\delta
x(t).\label{var-RDS}
\end{eqnarray}

In this paper, we first extend the concept of Hajnal diameter to
general matrices. A matrix with Hajnal diameter less than one has
the property of compressing the convex hull of
$\{x^{1},\cdots,x^{m}\}$. Consequently, for an infinite sequence
of time-varying Jacobian matrices, the average compression rate
can be used to verify synchronization. Since the Jacobian matrices
have identical row sums, the (skew) projection along the diagonal
synchronization direction can be used to define the projection
joint spectral radius, which equals the Hajnal diameter.
Furthermore, we show that the Hajnal diameter is equal to the
largest Lyapunov exponent along directions transverse to the
synchronization manifold; hence, it can also be used to determine
whether the coupled system (\ref{general}) can be synchronized.

Secondly, we apply these results to discuss the synchronization of
the CML with time-varying couplings. As we shall show, the Hajnal
diameter of infinite coupling stochastic matrices can be utilized
to measure the synchronizability of the coupling process. More
precisely,  the coupled system (\ref{TVCMLs}) synchronizes if the
sum of the logarithm of the Hajnal diameter and the largest
Lyapunov exponent of the uncoupled system is negative. Using the
equivalence of the Hajnal diameter, projection joint spectral
radius, and transverse Lyapunov exponents, we study some
particular examples for which the Hajnal diameter can be computed,
including static coupling, a finite coupling set, and a
multiplicative ergodic stochastic matrix process. We also present
numerical examples to illustrate our theoretical results.

The connection structure of the CML (\ref{CMLs}) naturally gives
rise to a graph, where each unit can be regarded as a vertex.
Hence, we associate the coupling matrix $G$ with a graph
$\Gamma=(V,E)$, with the vertex set $V=\{1,2,\dots,m\}$ and the
edge set $E=\{e_{ij}\}$, where there exists a directed edge from
vertex $j$ to vertex $i$ if and only if $G_{ij}>0$. The graphs we
consider here are assumed to be simple (that is, without loops and
multiple edges), but are allowed to be directed and weighted. That
is, we do not assume a symmetric coupling scheme.

We extend this idea to an infinite graph sequence $\{\Gamma(t)\}$.
That is, we regard a time-varying graph as a graph process
$\{\Gamma(t)\}_{t\in {\mathbb Z}^{+}}$. Define
$\Gamma(t)=[V,E(t)]$ where $V=\{1,2,\cdots,m\}$ denotes the vertex
set and $E(t)=\{e_{ij}(t)\}$ denotes the edge set of the graph at
time $t$. The time-varying coupling matrix $G(t)$ might then be
regarded as a function of the time-varying graph sequence, i.e.,
$G(t)=G(\Gamma(t))$. A basic problem that arises is, which kind of
sequence can ensure the synchrony of the coupled system for some
chaotic synchronized state $s(t+1)=f(s(t))$. As we shall show, the
property that the union of the $\Gamma(t)$ contains a spanning
tree is important for synchronizing chaotic maps. We prove that
under certain conditions, the coupling graph process can
synchronize some chaotic maps, if and only if there exists an
integer $T>0$ such that there exists at least one vertex $j$ from
which any other vertex can be accessible within a time interval of
length $T$.

This paper is organized as follows. In Section~2, we present some
definitions and lemmas on the Hajnal diameter, projection joint
spectral radius, projection Lyapunov exponents, and transverse
Lyapunov exponents for generalized Jacobian matrix sequences as
well as stochastic matrix sequences. In Section~3, we study the
synchronization of the generalized coupled discrete-time systems
with time-varying couplings (\ref{general}). In Section~4, we
discuss the synchronization of the CML with time-varying couplings
(\ref{TVCMLs})  and study the relation between synchronizability
and coupling graph process topologies. In addition, we present
some examples where synchronizability is analytically computable.
In Section~5, we present numerical examples to illustrate the
theoretical results, and conclude the paper in Section~6.

\section{Preliminaries}
In this section we present some definitions and lemmas on matrix
sequences. First, we extend the definitions of the Hajnal diameter
and the projection joint spectral radius, introduced in
\cite{Haj1,Haj2,Shen} for stochastic matrices, to generalized
time-varying matrix sequence. Furthermore, we extend Lyapunov
exponents and projection Lyapunov exponents to the general
time-varying case and discuss their relation. Secondly, we
specialize these definitions to stochastic matrix sequences and
introduce the relation between a stochastic matrix sequence and
graph topology.

\subsection{General definitions}
We  study the following generalized
time-varying linear system
\begin{eqnarray}
u(t+t_{0}+1)=L_{t+t_{0}}(\varrho^{(t)}(\phi))u(t+t_{0}),\label{linear}
\end{eqnarray}
where $\varrho^{(t)}$ is defined by a random dynamical system
$\{\Phi, {\mathcal B}, P,\varrho^{(t)}\}$, where $\Phi$ denotes
the state space, ${\mathcal B}$ the $\sigma$-algebra on $\Phi$,
$P$ the probability measure, $\varrho^{(t)}$ a semiflow. Studying
the linear system (\ref{linear}) comes from the variational system
of the coupled system (\ref{general}). For the variational system
(\ref{var-matrix}), $\varrho^{(t)}(\cdot)$ represents the
synchronized state flow $f^{(t)}(\cdot)$. And, if $L_{t}(\cdot)$
is independent of $t$, then the linear system (\ref{linear}) can
be rewritten as:
\begin{eqnarray}
u(t+1)=L(\varrho^{(t)}(\phi))u(t).\label{linear-RDS}
\end{eqnarray}
Thus, it can represent the variational system (\ref{var-RDS}) as a
special case, where $\varrho^{(t)}$ is the product flow
$(f^{(t)}(\cdot),\theta^{(t)}(\cdot))$. Hence, the linear system
(\ref{linear}) can unify the two cases of variational systems
(\ref{var-matrix},\ref{var-RDS}) of the coupled system
(\ref{general},\ref{RDS}).

For this purpose, we define a generalized matrix sequence map
${\mathcal L}$ from $\mathbb Z^{+}\times\Phi$ to $2^{\mathbb
R^{m\times m}}$,
\begin{eqnarray}
{\mathcal L}:{\mathbb Z}^{+}\times\Phi&\rightarrow& 2^{\mathbb{R}^{m\times m}}\nonumber\\
(t_{0},\phi)&\mapsto&
\{L_{t+t_{0}}(\varrho^{(t)}\phi)\}_{t\in{\mathbb Z}^{+}}.
\end{eqnarray}
where $2^{\mathbb R^{m\times m}}$ denotes the set containing all
subsets of $\mathbb R^{m\times m}$. In \cite{Haj1,Haj2}, the
concept of the Hajnal diameter was introduced to describe the
compression rate of a stochastic matrix. We extend it to general
matrices below.

\begin{definition}
\label{def2.1} For a matrix $L$ with row vectors
$g_{1},\cdots,g_{m}$ and a vector norm $\|\cdot\|$ in
$\mathbb{R}^{m}$, the Hajnal diameter of $L$ is defined by
\begin{eqnarray*}
{\rm diam}(L,\|\cdot\|)=\max\limits_{i,j}\|g_{i}-g_{j}\|.
\end{eqnarray*}
\end{definition}
We also introduce the Hajnal diameter for a matrix sequence map
${\mathcal L}$.
\begin{definition}
\label{def2.2} For a generalized matrix sequence map ${\mathcal
L}$, the Hajnal diameter of ${\mathcal L}$ at $\phi\in\Phi$ is
defined by
\begin{eqnarray*}
{\rm diam}({\mathcal
L},\phi)=\overline{\lim\limits_{t\rightarrow\infty}}\sup\limits_{t_{0}\ge
0}\left\{ {\rm
diam}(\prod\limits_{k=t_{0}}^{t_{0}+t-1}L_{k}(\varrho^{(k-t_{0})}\phi)\right\}^{\frac{1}{t}}.
\end{eqnarray*}
where $\prod$ denotes the left matrix product:
$\prod_{k=1}^{n}A_{k}=A_{n}\times A_{n-1}\times\cdots\times
A_{1}$.
\end{definition}

The Hajnal diameter for the infinite matrix
sequence map ${\mathcal L}$ does not depend on the choice of the
norm. In fact, all norms in a Euclidean space are equivalent and
any additional factor is eliminated by the power $1/t$ and the
limit as $t\rightarrow\infty$.

Let ${\mathcal H}\subset \mathbb{R}^{m\times m}$ be a class of
matrices having the property that all row sums are the same. Thus,
all matrices in ${\mathcal H}$ share the common eigenvector
$e_{0}=[1,1,\cdots,1]^{\top}$, where the corresponding eigenvalue
is the row sum of the matrix. Then, the projection joint spectral
radius can be defined for a generalized matrix sequence map
${\mathcal L}$, similar to introduced in \cite{Shen} as follows.

\begin{definition}\label{def2.3} Suppose ${\mathcal L}(t_{0},\phi)\subset{\mathcal H}$
for $t_{0}\in\mathbb Z^{+}$ and $\phi\in\Phi$. Let $\mathcal
E_{0}$ be the subspace spanned by the synchronization direction
$e_{0}=[1,1,\cdots,1]^{\top}$, and $P$ be any $(m-1)\times m$
matrix with  exact kernel  $\mathcal E_{0}$. We denote by
$\hat{L}\in\mathbb R^{(m-1)\times(m-1)}$ the (skew) projection of
matrix $L\in\mathcal H$ as the unique solution of
\begin{eqnarray}
PL=\hat{L} P.\label{proj}
\end{eqnarray}
The projection joint spectral radius of the generalized matrix
sequence map ${\mathcal L}$ is defined as
$$
\hat{\rho}({\mathcal
L},\phi)=\overline{\lim\limits_{t\rightarrow\infty}}\sup\limits_{t_{0}\ge
0}{\bigg\|\prod\limits_{k=t_{0}}^{t_{0}+t-1}\hat{L}_{k}(\varrho^{(k-t_{0})}\phi)\bigg\|}
^{\frac{1}{t}}.$$
\end{definition}

One can see that $\hat{\rho}({\mathcal L},\phi)$ is
independent of the choice of the matrix norm $\|\cdot\|$ induced
by vector norm. The following lemma shows that it is also
independent of the choice of the matrix $P$.

\begin{lemma}\label{lem2.4} Suppose ${\mathcal L}(t_{0},\phi)\subset {\mathcal H}$
for all $t_{0}\ge 0$ and $\phi\in\Phi$. Then
$$ \hat{\rho}({\mathcal
L},\phi)={\rm diam}({\mathcal L},\phi).$$
\end{lemma}
A proof is given in the Appendix.

The Lyapunov exponents are often used to study evolution of the
dynamics \cite{Jost,Ding}. Here, we extend the definitions of
Lyapunov exponents to general time-varying cases.

\begin{definition}\label{def2.5}
For the coupled system (\ref{general}), the Lyapunov exponent of
the matrix sequence map ${\mathcal L}$ initiated by $\phi\in\Phi$
in the direction $u\in \mathbb{R}^{m}$ is defined as
\begin{eqnarray}
\lambda(\mathcal{L},\phi,u)=\overline{\lim\limits_{t\rightarrow\infty}}\frac{1}{t}\sup\limits_{t_{0}\ge
0}
\log\big\|\prod\limits_{k=t_{0}}^{t+t_{0}-1}L_{k}(\varrho^{(k-t_{0})}\phi)u\big\|.
\label{genanralyap}
\end{eqnarray}
The projection along the synchronization direction $e_{0}$ can
also define a Lyapunov exponent, called the projection Lyapunov
exponent:
\begin{eqnarray}
\hat{\lambda}({\mathcal
L},\phi,v)=\overline{\lim\limits_{t\rightarrow\infty}}
\frac{1}{t}\sup\limits_{t_{0}\ge
0}\log\big\|\prod\limits_{k=t_{0}}^{t+t_{0}-1}\hat{L}_{k}(\varrho^{(k-t_{0})}\phi)v\big\|,
\end{eqnarray}
where $\hat{L}_{k}(\varrho^{k}\phi)$ is the projection of matrix
$L_{k}(\varrho^{k}\phi)$ as defined in Definition \ref{def2.3}.
\end{definition}

It can be seen that the definition of the generalized Lyapunov
exponent above satisfies the basic properties of Lyapunov
exponents\footnote{This kind of definition of characteristic
exponent is similar to the Bohl exponent used to study uniform
stability of time-varying systems in \cite{Bohl}.}. For more
details about generalized Lyapunov exponents, we refer to
\cite{Barr}.
\begin{lemma}\label{lem2.6}
Suppose ${\mathcal L}(t_{0},\phi)\subset{\mathcal H}$ for all
$\phi\in\Phi$ and $t_{0}\ge 0$. Then,
$$\sup\limits_{v\in \mathbb{R}^{m-1},v\ne 0}\hat{\lambda}({\mathcal
L},\phi,v)=\log\hat{\rho}({\mathcal L},\phi)=\log{\rm
diam}(\mathcal L,\phi).$$
\end{lemma}
A proof is given in the Appendix.

This lemma implies that the projection joint spectral radius gives
the largest Lyapunov exponent in directions transverse to the
synchronization direction $e_{0}$ of the matrix sequence map
${\mathcal L}$.

When the time dependence arises from being totally driven by some
random dynamical system, we can write the generalized matrix
sequence map $\mathcal L$ as $\mathcal
L(\phi)=\{L(\varrho^{(t)}\phi)\}_{t\in\mathbb Z^{+}}$ since it is
independent of $t_{0}$ and is just a map on $\Phi$. As introduced
in \cite{Col}, we have specific definitions for Lyapunov exponents
of the time-varying system (\ref{linear-RDS}) as follows.

For the linear system (\ref{linear-RDS}), the Lyapunov exponent of
the matrix sequence map ${\mathcal L}$ initiated by $\phi\in\Phi$
in the direction $u\in \mathbb{R}^{m}$ is defined as
\begin{eqnarray}
\lambda(\mathcal{L},\phi,u)=\overline{\lim\limits_{t\rightarrow\infty}}\frac{1}{t}
\log\bigg\|\prod\limits_{k=0}^{t-1}L(\varrho^{(k)}\phi)u\bigg\|.
\label{Lyap}
\end{eqnarray}

If ${\mathcal L}(\phi)\subset{\mathcal H}$ for all $\phi\in\Phi$,
then the Lyapunov exponent in the synchronization direction
$e_{0}$ is
\begin{eqnarray}
\lambda({\mathcal
L},\phi,e_{0})=\overline{\lim\limits_{t\rightarrow\infty}}\frac{1}{t}\log\sum\limits_{k=0}^{t-1}|c(k)|,
\end{eqnarray}
where $c(k)$ denotes the corresponding common row sum at each time
$k$. The projection along the synchronization direction $e_{0}$
can also define a Lyapunov exponent, called the projection
Lyapunov exponent:
\begin{eqnarray}
\hat{\lambda}({\mathcal
L},\phi,v)=\overline{\lim\limits_{t\rightarrow\infty}}
\frac{1}{t}\log\bigg\|\prod\limits_{k=0}^{t-1}\hat{L}_{k}(\varrho^{(k)}\phi)v\bigg\|,
\end{eqnarray}
where $\hat{L}(\varrho^{k}\omega)$ is the (skew) projection of
matrix $L(\varrho^{k}\omega)$. Also, the Hajnal diameter and
projection joint spectral radius become
\begin{eqnarray*}
{\rm diam}(\mathcal
L,\phi)=\overline{\lim\limits_{t\rightarrow\infty}}\bigg\{{\rm
diam}\big(\prod\limits_{k=0}^{t-1}L(\varrho^{(t)}\phi)\big)\bigg\}^{\frac{1}{t}},\quad
\hat{\rho}(\mathcal
L,\phi)=\overline{\lim\limits_{t\rightarrow\infty}}\bigg\|\prod\limits_{k=0}^{t-1}\hat{L}(\varrho^{(k)}\phi)\bigg\|
^{\frac{1}{t}}.
\end{eqnarray*}
According to Lemmas \ref{lem2.4} and \ref{lem2.6}, $\log{\rm
diam}(\mathcal L,\phi)=\log\hat{\rho}(\mathcal
L,\phi)=\sup_{v\in\mathbb R^{m-1},v\ne 0}\hat{\lambda}({\mathcal
L},\phi,v)$. Let $\lambda_{0}$ be the Lyapunov exponent along the
synchronization direction $e_{0}$ and $\lambda_{1}$,
$\lambda_{2}$, $\cdots$, $\lambda_{m-1}$ be the remaining Lyapunov
exponents for the initial condition $\phi$, counted with
multiplicities.

\begin{lemma}
\label{lem2.7} Suppose that ${\mathcal L}(\phi)\subset{\mathcal
H}$ is time-independent. Let the matrix
$D(t)=[D_{ij}(t)]_{i,j=1}^{m}$ denote the matrix
$L(\varrho^{(t)}\phi)$ and $c(t)$ denote the corresponding common
row sum of $D(t)$. If the following hold
\begin{enumerate}

\item
$\lim\limits_{t\rightarrow\infty}1/t\sum_{k=0}^{t-1}\log|c(k)|=\lambda_{0}$,

\item
$\overline{\lim\limits_{t\rightarrow\infty}}1/t\log^{+}|D_{ij}(t)|\le
0$, for all $i,j=\onetom$, where $\log^{+}(z)=\max\{\log z,0\}$,
\end{enumerate}
then
$$\log{\rm diam}(\mathcal L,\phi)=\log\hat{\rho}({\mathcal L},\phi)=\sup\limits_{i\ge
1}\lambda_{i}.$$
\end{lemma}
A proof is given in the Appendix.

Using the concept of Hajnal diameter, we can define (uniform)
synchronization of the non-autonomous system (\ref{general}) as
follows:
\begin{definition}
\label{def2.8} The coupled system (\ref{general}) is said to be
(uniformly locally completely) synchronized if there exists
$\eta>0$ such that for any $\epsilon>0$, there exists $T>0$ such
that the inequality
\begin{eqnarray}
{\rm diam}\big([x^{1}(t),x^{2}(t),\cdots,x^{m}(t)]^{\top}\big)\le
\epsilon
\end{eqnarray}
holds for all $t>t_{0}+T$, $t_{0}\ge 0$ and $x^{i}(t_{0})$,
$i=1,2,\cdots,m$ in the $\eta$ neighborhood of $s(t_{0})$ of a
synchronized state $s(t)$.
\end{definition}

\subsection{Stochastic matrix sequences}
The above definitions can also be used to deal with stochastic
matrix sequences.
\begin{definition}\label{def2.9}
A matrix $G\in \mathbb{R}^{m\times m}$ is said to be a stochastic
matrix if its elements are nonnegative and each row sum is $1$.
\end{definition}

We here consider the general time-varying case without the
assumption of an underlying random dynamical system and write a
stochastic matrix sequence as ${\mathcal G}=\{G(t)\}_{t\in{\mathbb
Z}^{+}}$. The case that the time variation is driven by some
dynamical system can be regarded as a special one.

\begin{definition}\label{def2.10}
The Hajnal diameter of ${\mathcal G}$ is defined as
\begin{eqnarray}
{\rm diam}({\mathcal
G})=\overline{\lim\limits_{t\rightarrow\infty}}\sup\limits_{t_{0}\ge
0}\bigg({\rm
diam}\prod\limits_{k=t_{0}}^{t_{0}+t-1}G(k)\bigg)^{\frac{1}{t}}
\end{eqnarray}
and the projection joint spectral radius for ${\mathcal G}$ is
\begin{eqnarray}
\hat{\rho}({\mathcal
G})=\overline{\lim\limits_{t\rightarrow\infty}}\sup\limits_{t_{0}\ge
0}\bigg\|\prod\limits_{k=t_{0}}^{t_{0}+t-1}\hat{G}(k)\bigg\|^{\frac{1}{t}}
\end{eqnarray}
where $\hat{G}(t)$ is the projection of $G(t)$, as in Definition \ref{def2.3}.
\end{definition}

Then, from Lemma \ref{lem2.4}, we have

\begin{lemma}\label{lem2.11}
$ {\rm diam}({\mathcal G})=\hat{\rho}({\mathcal G}).$
\end{lemma}

To estimate the Hajnal diameter of a product of stochastic
matrices, we use the concept of scrambling introduced in
\cite{Shen}.

\begin{definition}\label{def2.12}
A stochastic matrix $G=[G_{ij}]_{i,j=1}^{m}\in \mathbb{R}^{m\times
m}$ is said to be scrambling if for any $i,j$, there exists an
index $k$ such that $G_{ik}\ne 0$ and $G_{jk}\ne 0$.
\end{definition}

For $g_{i}=[g_{i,1},\cdots,g_{i,m}]\in \mathbb{R}^{m}$ and
$g_{j}=[g_{j,1},\cdots,g_{j,m}]\in \mathbb{R}^{m}$, define
\begin{eqnarray*}
g_{i}\wedge
g_{j}=\big[\min(g_{i,1},g_{j,1}),\cdots,\min(g_{i,m},g_{j,m})\big].
\end{eqnarray*}
We use the following quantity introduced in \cite{Haj1,Haj2} to
measure scramblingness,
\begin{eqnarray*}
\eta(G)=\min\limits_{i,j}\|g_{i}\wedge g_{j}\|_{1},
\end{eqnarray*}
where, $\|\cdot\|_{1}$ is the norm given by
$\|x\|_{1}=\sum_{i=1}^{m}|x_{i}|$ for
$x=[x_{1},\cdots,x_{m}]\in \mathbb{R}^{m}$. It is clear that
$0\le\eta(G)\le 1$, and that $\eta(G)>0$ if and only if $G$ is
scrambling. Thus, the well-known Hajnal inequality has the
following generalized form.

\begin{lemma}\label{def2.13}~(Generalized Hajnal inequality, Theorem 6 in \cite{Shen}.)
For any vector norm in $\mathbb{R}^{m}$ and any two stochastic
matrices $G$ and $H$,
\begin{eqnarray}
{\rm diam}(GH)\le (1-\eta(G)){\rm diam}(H).
\end{eqnarray}
\end{lemma}

The concepts of projection joint spectral radius and Hajnal
diameter are linked to the ergodicity of stochastic matrix
sequences. We can extend the ergodicity for a matrix set
\cite{Shen,Wolf} to a matrix sequence as follows:

\begin{definition}\label{def2.14}  (Ergodicity, Definition 1 in \cite{Mor}.)
A stochastic matrix sequence $\Sigma=\{G(t)\}_{t\in{\mathbb
Z}^{+}}$ is said to be ergodic if for any $t_{0}$ and
$\epsilon>0$, there exists $T>0$ such that for any $t>T$ and some
norm $\|\cdot\|$,
\begin{eqnarray}
{\rm diam} \left( \prod\limits_{s=t_{0}}^{t_{0}+t-1} G(s) \right)\le
\epsilon.\label{erg}
\end{eqnarray}
Moreover, if for any $\epsilon>0$, there exists $T>0$ such that
inequality (\ref{erg}) holds for all $t\ge T$ and $t_{0}\ge 0$,
${\mathcal G}$ is said to be uniformly ergodic.
\end{definition}

A stochastic matrix $G=[G_{ij}]_{i,j=1}^{m}$ can be associated
with a graph $\Gamma=[V,E]$, where $V=\{1,2,\cdots,m\}$ denotes
the vertex set and $E=\{e_{ij}\}$ the edge set, in the sense that
there exists an edge from vertex $j$ to $i$ if and only if
$G_{ij}>0$. Let $\Gamma_{1}=[V,E_{1}]$ and $\Gamma=[V, E_{2}]$ be
two simple graphs with the same vertex set. We also define the
union $\Gamma_{1}\bigcup\Gamma_{2}=[V,E_{1}\bigcup E_{2}]$
(merging multiple edges). It can be seen that for two stochastic
matrices $G_{1}$ and $G_{2}$ with the same dimension and positive
diagonal elements, the edge set of $\Gamma_{1}\bigcup\Gamma_{2}$
is contained in that of the corresponding graph of the product
matrix $G_{1} G_{2}$. In this way, we can define the union of the
graph sequence $\{\Gamma(t)\}_{t\in\mathbb Z^{+}}$ across the time
interval $[t_{1},t_{2}]$ by
$\bigcup_{k=t_{1}}^{t_{2}}\Gamma(k)=[V,\bigcup_{k=t_{1}}^{t_{2}}E(k)]$.
The following concepts for graphs can be found, e.g., in
\cite{God}.

\begin{definition}\label{def2.15}
A graph $\Gamma$ is said to have a spanning tree if there exists a
vertex, called the root, such that for each other vertex $j$ there
exists at least one directed path from the root to vertex $j$.
\end{definition}

It follows that $\{\Gamma(t)\}_{t\in\mathbb Z^{+}}$ has a spanning
tree across the time interval $[t_{1},t_{2}]$ if the union of
$\{\Gamma(t)\}_{t\in\mathbb Z^{+}}$ across $[t_{1},t_{2}]$ has a
spanning tree. This is equivalent to the existence of a vertex
from which all other vertices can be accessible across
$[t_{1},t_{2}]$.

\begin{definition}\label{def2.16}
A graph $\Gamma$ is said to be scrambling if for any different
vertices $i$ and $j$, there exists a vertex $k$ such that there
exist edges from $k$ to $i$  and from $k$ to $j$.
\end{definition}

It follows that a stochastic matrix $G$ is scrambling if and only
if the corresponding graph $\Gamma$ is scrambling.

\begin{lemma}\label{lem2.17}~(See Lemma 4 in \cite{Wolf}.)
Let $G(1),G(2),\cdots,G(m-1)$ be stochastic matrices with positive
diagonal elements, where each of the corresponding graphs
$\Gamma(1)$, $\Gamma(2)$, $\cdots$, $\Gamma(m-1)$ have spanning
trees. Then $\prod_{k=1}^{m-1}G(k)$ is scrambling.
\end{lemma}

Suppose now that the stochastic matrix sequence ${\mathcal G}$ is
driven by some metric dynamical system $\mathcal
Y=\{\Omega,{\mathcal F},P,\theta^{(t)}\}$. We write ${\mathcal G}$
as $\{G(t)=G(\theta^{(t)}\omega)\}_{t\in\mathbb Z^{+}}$, where
$\omega\in\Omega$. Then, as stated in Section 2.1, we can define
the Lyapunov exponents.
\begin{definition}\label{def2.18}
The Lyapunov exponent of the stochastic matrix sequence $\mathcal
G$ is defined as
\begin{eqnarray*}
\sigma(\mathcal
G,\omega,u)=\overline{\lim\limits_{t\rightarrow\infty}}
\frac{1}{t}\log\big\|\prod\limits_{k=0}^{t-1}G(\theta^{t}\omega)u\big\|.
\end{eqnarray*}
The projection Lyapunov exponents is defined as
\begin{eqnarray*}
\hat{\sigma}(\mathcal
G,\omega,u)=\overline{\lim\limits_{t\rightarrow\infty}}
\frac{1}{t}\log\big\|\prod\limits_{k=0}^{t-1}\hat{G}(\theta^{t}\omega)u\big\|,
\end{eqnarray*}
where $\hat{G}(\cdot)$ is the projection of $G(\cdot)$ as defined
in Definition \ref{def2.3}.
\end{definition}

For a given $\omega\in\Omega$, one can see that ${\rm
diam}({\mathcal G})$ and $\hat{\rho}({\mathcal G})$ both equal the
largest Lyapunov exponent of ${\mathcal G}$ in directions
transverse to the synchronization direction under several mild
conditions.

In closing this section, we list some notations to be used in the
remainder of the paper. The matrix $\hat{L}$ denotes the (skew)
projection of the matrix $L$ along the vector $e$ introduced in
Definition \ref{def2.3}, and $\hat{\mathcal L}$ is the (skew)
projection of the matrix sequence map ${\mathcal L}$ along $e$. For
$x=(x^{1},\cdots,x^{m})^{\top}\in \mathbb{R}^{m}$, the average
$\frac{1}{m}\sum_{i=1}^{m}x^{i}$  of $x$ is denoted by $\bar{x}$.
The notation $\|\cdot\|$ denotes some vector norm in the linear
space $\mathbb{R}^{m}$, and also the matrix norm in
$\mathbb{R}^{m\times m}$ induced by this vector norm.
$f^{(t)}(s_{0})$ denotes the $t$-iteration of the map $f$ with
initial condition $s_{0}$. We let $x(t,t_{0},x_{0})$ be the solution
of the coupled system (\ref{general}) with initial condition
$x(t_{0})=x_{0}$, which we sometimes abbreviate as $x(t)$.

\section{Generalized synchronization analysis}
For the variational system (\ref{var-matrix}), similar to the
Subsection 2.1, we denote by $\mathcal D$ the Jacobian sequence
map in the generalized sense, i.e., $\mathcal D$ is a map from
$\mathbb Z^{+}\times\mathbb R$ to $2^{\mathbb R^{m\times m}}$:
$\mathcal
D(t_{0},s_{0})=\{D_{t+t_{0}}(f^{(t)}(s_{0}))\}_{t\in{\mathbb
Z}^{+}}\subset\mathcal H$ for all $t_{0}\in\mathbb Z^{+}$ and
$s_{0}\in A$. Furthermore, letting
$$B(t,t_{0})=\prod\limits_{k=t_{0}}^{t+t_{0}-1}D_{k}(f^{(k-t_{0})}(s_{0})),$$
we can rewrite the variational system (\ref{var-matrix}) as
follows:
\begin{eqnarray}
\delta x(t+t_{0})&=&D_{t+t_{0}-1}(f^{(t-1)}(s_{0}))\delta
x(t+t_{0}-1)=B(t,t_{0})\delta x(t_{0}).\label{var}
\end{eqnarray}

From Definitions \ref{def2.2} and \ref{def2.3}, we have
\begin{eqnarray*}
{\rm diam}({\mathcal
D},s_{0})&=&\overline{\lim\limits_{t\rightarrow\infty}}
\sup\limits_{t_{0}\ge
0}\bigg\{{\rm diam}\big(\prod\limits_{k=t_{0}}^{t_{0}+t-1}D_{k}(f^{(k-t_{0})}(s_{0}))\big)\bigg\}^{\frac{1}{t}},\\
\hat{\rho}({\mathcal
D},s_{0})&=&\overline{\lim\limits_{t\rightarrow\infty}}\sup\limits_{t_{0}\ge
0}\bigg\|\prod\limits_{k=t_{0}}^{t_{0}+t-1}\hat{D}_{k}(f^{(k-t_{0})}(s_{0}))\bigg\|^{\frac{1}{t}}.
\end{eqnarray*}
We will also refer to the following hypothesis.

{\em$\bf H_{3}$.
\begin{eqnarray}
\sup\limits_{s_{0}\in{A}}{\rm diam}({\mathcal
D},s_{0})<1.\label{thm1}
\end{eqnarray}}

\begin{theorem}\label{thm3.1}
If hypotheses $\bf H_{1}$--$\bf H_{3}$ hold, then the compact set
$A^{m}\bigcap\mathcal S$ is a uniformly asymptotically stable
attractor of the coupled system (\ref{general}) in
$\mathbb{R}^{m}$, i.e., the coupled system (\ref{general}) is
uniformly locally completely synchronized.
\end{theorem}
{\em Proof.} Let
\begin{eqnarray*}
{\rm diam}({\mathcal D},t_{0},t,s_{0})
={\rm diam}\bigg(\prod\limits_{k=t_{0}}^{t_{0}+t-1}D_{k}\big(f^{(k-t_{0})}(s_{0})\big)\bigg),\\
{\rm diam}({\mathcal D},t,s_{0})=\sup\limits_{t_{0}\ge 0}
\bigg\{{\rm
diam}\big(\prod\limits_{k=t_{0}}^{t_{0}+t-1}D_{k}(f^{(k-t_{0})}(s_{0}))\big)\bigg\}.
\end{eqnarray*}
According to $\bf H_{3}$, letting $1>d>\sup_{s_{0}\in{A}}{\rm
diam}({\mathcal D},s_{0})$ and $n_{0}$ satisfy
$d^{n_{0}}<\frac{1}{3}$, for any $s_{0}\in A$, there exists
$n(s_{0})\ge n_{0}$ such that ${\rm diam}({\mathcal D},t,s_{0})<d$
holds for all $t\ge n(s_{0})$. By equicontinuity ($\bf H_{1}$) and
compactness ($\bf H_{2}$), there must exist a finite integer set
$\mathcal V=\{n_{1},n_{2},\cdots,n_{v}\}$ satisfying $n_{i}\ge
n_{0}$ for all $i=1,2,\cdots,v$ and a neighborhood $U$ of $A$ such
that for any $s_{0}\in U$, there exists $n_{j}\in\mathcal V$ such
that ${\rm
diam}\big(\prod_{k=t_{0}}^{t_{0}+n_{j}-1}D_{k}(f^{(k-t_{0})}(s_{0}))\big)<d^{n_{j}}<\frac{1}{3}$
holds for all $t_{0}\ge 0$.

By the hypothesis $\bf H_{2}$, there exists a compact neighborhood
$W$ of $A$ such that $U\supset W\supset A$, $f(W)\subset W$, and
$\bigcap\limits_{n\ge 0}f^{(n)}(W)=A$ \cite{Mil}. Let
\begin{eqnarray*}
a=\min\limits_{n\in\mathcal V}d_{H}(f^{(n)}(W),W)>0,
\end{eqnarray*}
where $d_{H}(\cdot,\cdot)$ denotes the Hausdorff metric in
$\mathbb R$ . Then, define a compact set
\begin{eqnarray*}
W_{\alpha}=\bigg\{x=(x^{1},\cdots,x^{m})\in
\mathbb{R}^{m}:~\max_{1\le i\le
m}|x^{i}-\bar{x}|\le\alpha~and~\bar{x}\in W\bigg\}.
\end{eqnarray*}
By the mean value theorem, we have
\begin{eqnarray*}
f^{i}_{k}(x^{1}(k),\cdots,x^{m}(k))-f(s(k),\cdots,s(k))=\sum\limits_{j=1}^{m}\frac{\partial
f_{k}^{i}}{\partial x^{j}}(\xi^{ij}_{k}),
\end{eqnarray*}
where $\xi^{ij}_{k}$ belongs to the closed interval induced by the
two ends $x^{i}(k)$ and $s(k)$. Denote by $D_{k}(\xi_{k})$ the
matrix $[\partial f_{k}^{i}(\xi^{ij}_{k})/\partial
x^{j}]_{i,j=1}^{m}$.

Let $\alpha>0$ be sufficiently small so that for each $x_{0}\in
W_{\alpha}$ with $s(t_{0})=\bar{x}_{0}$ and $x(t_{0})=x_{0}$,
there exists $t_{1}\in\mathcal V$ such that
\begin{eqnarray*}
&&|x^{i}(t_{1},t_{0},x_{0})-f^{(t_{1}-t_{0})}(\bar{x}_{0})|\le\frac{a}{2}\\
&&{\rm
diam}\bigg(\prod\limits_{k=t_{0}}^{t_{0}+t_{1}-1}D_{k}(\xi_{k})\bigg)<\frac{1}{2}
\end{eqnarray*}
holds for all $t_{0}\ge 0$. Then, for any $x_{0}\in W_{\alpha}$,
$\bar{x}_{0}\in W$, we have
\begin{eqnarray*}
\delta
x(t_{1}+t_{0})=\prod\limits_{k=t_{0}}^{t_{1}+t_{0}-1}D_{k}(\xi(k))\delta
x_{0}=\tilde{B}(t_{1},t_{0})\delta x_{0},
\end{eqnarray*}
where
$\tilde{B}(t_{1},t_{0})=\prod_{k=t_{0}}^{t_{1}+t_{0}-1}D_{k}(\xi(k))$.
Then,
\begin{eqnarray*}
|\delta x^{i}(t_{1}+t_{0})-\delta
x^{j}(t_{1}+t_{0})|&\le&\sum\limits_{k=1}^{m}|\tilde{B}_{ik}(t_{1},t_{0})-\tilde{B}_{jk}(t_{1},t_{0})||\delta
x^{j}_{0}|\\
&\le&{\rm diam}(\tilde{B}(t_{1},t_{0}))\max_{1\le i\le
m}|x_{0}^{i}-\bar{x}_{0}|.
\end{eqnarray*}
Thus, we conclude that
\begin{eqnarray*}
\max_{1\le i,j\le m}|x^{i}(t_{1}+t_{0})-x^{j}(t_{1}+t_{0})|\le
\frac{1}{2} \max_{1\le i,j\le m}|x^{i}_{0}-x^{j}_{0}|.
\end{eqnarray*}
By the definition of $W_{\alpha}$, we see that $x(t_{1}+t_{0})\in
W_{\alpha/2}$. With initial time $t_{0}+t_{1}$, we can continue
this phase and afterwards obtain
\begin{eqnarray*}
\lim\limits_{t\rightarrow\infty}|x^{i}(t)-x^{j}(t)|=0,~i,j=\onetom,
\end{eqnarray*}
uniformly with respect to $t_{0}\in\mathbb Z^{+}$ and $x_{0}\in
W_{\alpha}$. Therefore, the coupled system (\ref{general}) is
uniformly synchronized. Furthermore, we obtain that
$A^{m}\bigcap\mathcal S$ is a uniformly asymptotically stable
attractor for the coupled system (\ref{general}) and the convergence
rate can be estimated by $O(\{\sup_{s_{0}\in A}{\rm diam}({\mathcal
D},s_{0})\}^{t})$ since $d$ is chosen arbitrarily greater than
$\sup_{s_{0}\in A}{\rm diam}(\mathcal D,s_{0})$. The theorem is
proved.

{\sc Remark 1.} The idea of the above proof comes from that
of Theorem 2.12 in \cite{Ash}, with a modification for the
time-varying case. In Theorem 2.12 in \cite{Ash}, the authors used
normal Lyapunov exponents to prove asymptotical stability of the
original autonomous  system for the case when it is asymptotically
stable in an invariant manifold. In this paper, we directly use
the Hajnal diameter of the left product of the infinite Jacobian
matrix sequence map to measure the transverse differences of the
collections of spatial states. Furthermore, we consider a
non-autonomous system here due to time-varying couplings.

Following Lemma \ref{lem2.4} gives
\begin{corollary}\label{col3.2}
If $\sup_{s_{0}\in A}\hat{\rho}({\mathcal D},s_{0})<1$, then the
coupled system (\ref{general}) is uniformly synchronized.
\end{corollary}

Consider the special case that the coupled system (\ref{RDS}) is a
RDS on a MDS $\mathcal Y=\{\Omega,{\mathcal F},P,\theta^{(t)}\}$.
We can write this coupled system (\ref{RDS}) as a product
dynamical system $\{A\times\Omega,{\bf F }, {\bf
P},\Theta^{(t)}\}$, where ${\bf F}$ is the product
$\sigma$-algebra on $A\times\Omega$, ${\bf P}$ denotes the
probability measure, and
$\Theta^{(t)}(s_{0},\omega)=(\theta^{(t)}\omega,f^{(t)}(s_{0}))$.
Let $D(f^{(t)}(s_{0}),\theta^{(t)}\omega)$ denote the Jacobian
matrix at time $t$. By Definition \ref{def2.5}, the Lyapunov
exponents for the coupled system (\ref{RDS}) can be written as
follows:
\begin{eqnarray*}
\lambda(u,s_{0},\omega)=\overline{\lim\limits_{t\rightarrow\infty}}\frac{1}{t}\log\bigg\|\prod\limits_{k=0}^{t-1}D(
f^{(k)}(s_{0}),\theta^{(k)}\omega)u\bigg\|.
\end{eqnarray*}
It can be seen that the Lyapunov exponent along the diagonal
synchronization direction $e_{0}$ is
\begin{eqnarray*}
\lambda(e_{0},s_{0},\omega)=\overline{\lim\limits_{t\rightarrow\infty}}\frac{1}{t}\sum\limits_{k=0}^{t-1}
\log|c(k)|,
\end{eqnarray*}
where $c(k)$ is the common row sum of $D(
f^{(k)}(s_{0}),\theta^{(k)}\omega)$. Let
$\lambda_{0}=\lambda(e_{0},s_{0},\omega)$, $\lambda_{1}$,
$\cdots$, $\lambda_{m-1}$ be the Lyapunov exponents (counting
multiplicity) of the dynamical system ${\mathcal L}$ with the
initial condition $(s_{0},\omega)$. From Lemma \ref{lem2.7}, we
conclude that $\sup_{i\ge
1}\lambda_{i}=\log\hat{\rho}(F,s_{0},\omega)=\log {\rm
diam}(F,s_{0},\omega)$. If the probability ${\bf P}$ is ergodic,
then the Lyapunov exponents exist for almost all $s_{0}\in A$ and
$\omega\in\Omega$, and furthermore they are independent of
$(s_{0},\omega)$.

\begin{corollary}\label{col3.3} Suppose that hypotheses $\bf H_{1}$-$\bf H_{2}$
and the assumptions in Lemma \ref{lem2.7} hold. Suppose further
that $A\times\Omega$ is compact in the weak topology defined in
this RDS, the semiflow $\Theta^{(t)}$ is continuous, the Jacobian
matrix $D(\cdot,\cdot)$ is non-singular and continuous on
$A\times\Omega$, and
$$\sup\limits_{{\bf P}\in \rm{Erg}_{\Theta}(A\times\Omega)}\sup\limits_{i\ge 1}\lambda_{i}<0,$$
where $\rm{Erg}_{\Theta}(A\times\Omega)$ denotes the ergodic
probability measure set supported in $\{A\times\Omega,{\bf F
},\Theta^{(t)}\}$. Then the coupled system (\ref{RDS}) is
uniformly locally completely synchronized.
\end{corollary}

{\em Proof.} By Theorem 2.8 in \cite{Ash}, we have
\begin{eqnarray*}
\sup\limits_{{\bf P}\in
Erg_{\Theta}(A\times\Omega)}\lambda_{\max}(\hat{\mathcal D},{\bf
P})=\sup\limits_{\|u\|=1,(s_{0},\omega)\in
A\times\Omega}\overline{\lim\limits_{t\rightarrow\infty}}
\frac{1}{t}\log\bigg\|\prod\limits_{k=0}^{t-1}\hat{D}(f^{(k)}(s_{0}),\theta^{(k)}\omega)u\bigg\|,
\end{eqnarray*}
where $\hat{\mathcal D}$ is the projection of the intrinsic matrix
sequence map ${\mathcal D}$ and $\lambda_{\max}(\hat{\mathcal
D},{\bf P})$ denotes the largest Lyapunov exponent of
$\hat{\mathcal D}$ according to the ergodic probability ${\bf P}$
(the value for all almost $(s_{0},\omega)$ according to ${\bf
P}$). From Lemmas \ref{lem2.4}, \ref{lem2.6} and \ref{lem2.7}, it
follows
\begin{eqnarray*}
\sup\limits_{{\bf P}\in
\rm{Erg}_{\Theta}(A\times\Omega)}\sup\limits_{i\ge
1}\lambda_{i}&=&\sup\limits_{{\bf P}\in
\rm{Erg}_{\Theta}(A\times\Omega)}\lambda_{\max}(\hat{\mathcal
D},{\bf P})=\sup\limits_{(s_{0},\omega)\in
A\times\Omega}\lambda_{\max}(\hat{\mathcal
D},s_{0},\omega)\\
&=&\sup\limits_{(s_{0},\omega)\in
A\times\Omega}\log\hat{\rho}(\mathcal
D,s_{0},\omega)=\sup\limits_{(s_{0},\omega)\in
A\times\Omega}\log{\rm diam}(\mathcal D,s_{0},\omega).
\end{eqnarray*} The corollary is proved as a direct consequence from Theorem
\ref{thm3.1}.

{\sc Remark 2.} If $\lambda_{0}$ is the largest Lyapunov exponent,
then $V=\{u:~\lambda(u)<\lambda_{0}\}$ constructs a subspace of
$\mathbb R^{m}$ which is transverse to the synchronization
direction $e_{0}$. Corollary \ref{col3.3} implies that if all
Lyapunov exponents in the transverse directions are negative, then
the coupled system (\ref{general}) is synchronized. Otherwise, if
$\lambda_{0}$ is not the largest Lyapunov exponent, then
$\sup_{i\ge 1}\lambda_{i}<0$ implies that the largest exponent is
negative, which means that the synchronized solution $s(t)$ is
itself asymptotically stable through the evolution (\ref{RDS}).

{\sc Remark 3.} From Lemma \ref{lem2.7}, it can also be seen that
when computing $\rho({\mathcal D})$, it is sufficient to compute the
largest Lyapunov exponent of $\hat{\mathcal D}$.  In \cite{Ash}, the
authors proved for an autonomous dynamical system that if all
Lyapunov exponent of the normal directions, namely, the Lyapunov
exponents for $\hat{\mathcal D}$ are negative, then the attractor in
the invariant submanifold is an attractor in $\mathbb R^{m}$ (or a
more general manifold). In this paper, we extend the proof theorem
2.12 in \cite{Ash} to the general time-varying coupled system
(\ref{general})
 by discussing the relation between the Hajnal
diameter and transverse Lyapunov exponents. In the following
sections, we continue the synchronization analysis for
non-autonomous dynamical systems.

\section{Synchronization analysis of coupled map lattices with time-varying topologies}

Consider the following coupled system with time-varying
topologies:
\begin{eqnarray}
x^{i}(t+1)=\sum\limits_{j=1}^{m}G_{ij}(t)f(x^{j}(t)),i=\onetom,~t\in\mathbb
Z^{+},\label{lcmls}
\end{eqnarray}
where $f(\cdot): \mathbb{R}\rightarrow \mathbb{R}$ is $C^{1}$
continuous and $G(t)=[G_{ij}(t)]_{i,j=1}^{m}$ is a stochastic
matrix. In matrix form,
\begin{eqnarray}
x(t+1)=G(t)F(x(t)).\label{matrixlcmls}
\end{eqnarray}
Since the coupling matrix $G(t)$ is a stochastic matrix, the
diagonal synchronization manifold is invariant and we have the
uncoupled (or synchronized) state as:
\begin{eqnarray}
s(t+1)=f(s(t)).\label{uncoupled}
\end{eqnarray}
We suppose that for the synchronized state (\ref{uncoupled}),
there exists an asymptotically stable attractor $A$ with the
(maximum) Lyapunov exponent
\begin{eqnarray*}
\mu=\sup\limits_{s_{0}\in
A}\overline{\lim\limits_{t\rightarrow\infty}}\frac{1}{t}\sum\limits_{k=0}^{t-1}\log|f'(s(k))|.
\end{eqnarray*}
The system (\ref{lcmls}) is a special form of (\ref{general})
satisfying the equicontinuous condition $\bf H_{1}$. Linearizing
the system (\ref{lcmls}) about the synchronized state yields the
variational equation
\begin{eqnarray*}
\delta x^{i}(t+1)=\sum\limits_{j=1}^{m}G_{ij}(t)f'(s(t))\delta
x^{i}(t),~i=\onetom,
\end{eqnarray*}
and
\begin{eqnarray*}
{\rm
diam}\bigg(\prod\limits_{k=t_{0}}^{t_{0}+t-1}G(k)f'(f^{(k-t_{0})}(s_{0}))\bigg)
={\rm
diam}\bigg(\prod\limits_{k=t_{0}}^{t_{0}+t-1}G(k)\bigg)\bigg|\prod\limits_{l=0}^{t}f'(f^{(l)}(s_{0}))\bigg|.
\end{eqnarray*}
Denote the stochastic matrix sequence $\{G(t)\}_{t\in\mathbb
Z^{+}}$ by $\mathcal G$. Thus, the Hajnal diameter of the
variational system is ${\rm diam}({\mathcal G})e^{\mu}$. Using
Theorem \ref{thm3.1}, we have the following result.

\begin{theorem}\label{thm4.1}
Suppose that the uncoupled system $s(t+1)=f(s(t))$ satisfies
hypothesis $\bf H_{2}$ with Lyapunov exponent $\mu$. Let
${\mathcal G}=\{G(t)\}_{t\in{\mathbb Z}^{+}}$. If
\begin{eqnarray}
{\rm diam}({\mathcal G})e^{\mu}<1,\label{thm2}
\end{eqnarray}
then the coupled system (\ref{lcmls}) is synchronized.
\end{theorem}

From Theorem \ref{thm4.1}, one can see that the quantity ${\rm
diam}({\mathcal G})$ as well as other equivalent quantities such
as the projection joint spectral radius and the Lyapunov exponent,
can be used to measure the synchronizability of the time-varying
coupling, i.e., the coupling stochastic matrix sequence ${\mathcal
G}$. A smaller value of ${\rm diam}({\mathcal G})$ implies a
better synchronizability of the time-varying coupling topology. If
the uncoupled system (\ref{uncoupled}) is chaotic, i.e.~$\mu>0$,
then the necessary condition for synchronization condition
(\ref{thm2}) is ${\rm diam}({\mathcal G})<1$. So, it is important
to investigate under what conditions ${\rm diam}({\mathcal G})<1$
holds.

Suppose that the stochastic matrix set ${\mathcal M}$ satisfies the
following hypotheses:

{\em $\bf H_{4}$. ${\mathcal M}$ is compact and there exists $r>0$
such that for any $G=[G_{ij}]_{i,j=1}^{m}\in{\mathcal M}$,
$G_{ij}>0$ implies $G_{ij}\ge r$ and all diagonal elements
$G_{ii}>r$, $i=\onetom$.}

We denote the graph sequence corresponding to the stochastic
matrix sequence ${\mathcal G}$ by
${\mathbf\Gamma}=\{\Gamma(t)\}_{t\in{\mathbb Z}^{+}}$.
Then we have the following result.

\begin{theorem}\label{thm4.2}
Suppose that the stochastic matrix sequence ${\mathcal G}\subset
{\mathcal M}$  satisfies hypothesis $\bf H_{4}$. Then, the
following statements are equivalent:
\begin{enumerate}
\item ${\rm diam}({\mathcal G})<1$;

\item there exists $T>0$ such that for any $t_{0}$, the graph
$\bigcup_{k=t_{0}}^{t_{0}+T}\Gamma(k)$ has a spanning tree;

\item the stochastic matrix sequence ${\mathcal G}$ is uniformly
ergodic.
\end{enumerate}
\end{theorem}

{\em Proof.} We first show $(3)\Rightarrow(2)$ by reduction to
absurdity. Let $B(t_{0},t)=\prod_{k=t_{0}}^{t_{0}+t-1}G(k)$. Since
${\mathcal G}$ is uniformly ergodic, there must exist $T>0$ such
that ${\rm diam}(B(t_{0},T))<1/2$ holds for any $t_{0}\ge 0$. So,
$v=\prod_{k=t_{0}}^{t_{0}+T-1}G(k)u$ satisfies:
\begin{eqnarray}
\max_{1\le i,j\le m}|v_{i}-v_{j}|\le {\rm
diam}(B(t_{0},T))\|u\|_{\infty}\le\frac{1}{2}\|u\|_{\infty}.\label{ineq}
\end{eqnarray}
If the second condition does not hold, then there exists $t_{T}$
such that the union $\bigcup_{k=t_{T}}^{t_{T}+T-1}\Gamma(k)$ does
not have a spanning tree. That is, there exist two vertices
$v_{1}$ and $v_{2}$ such that for any vertex $z$, there is either
no directed path from $z$ to $v_{1}$ or no directed path from $z$
to $v_{2}$. Let $U_{1}$ ($U_{2}$) be the vertex set which can
reach $v_{1}$ ($v_{2}$, respectively) across $[t_{T},t_{T}+T-1]$.
This implies that $U_{1}$ and $U_{2}$ are disjoint across
$[t_{T},t_{T}+T-1]$ and no edge starts outside of $U_{1}$
($U_{2}$) and ends in $U_{1}$ ($U_{2}$) . Furthermore, considering
the Frobenius form of $G(t)$, one can see that the elements in the
corresponding rows of $U_{1}$ ($U_{2}$) with columns associated
with outside of $U_{1}$ ($U_{2}$) are all zeros. Let
\begin{eqnarray*}
u_{i}=\left\{\begin{array}{ll} 1&i\in U_{1},\\
0&i\in U_{2},\\
\text{any~value~in~}(0,1),& \text{otherwise}.\end{array}\right.
\end{eqnarray*}
We have
\begin{eqnarray*}
v_{i}=\left\{\begin{array}{ll}1&i\in U_{1},\\
0&i\in U_{2},\\
\in[0,1],&\text{otherwise}.\end{array}\right.
\end{eqnarray*}
This implies that $\max_{1\le i,j\le m}|v_{i}-v_{j}|\ge
1=\|u\|_{\infty}$, which contradicts with (\ref{ineq}). Therefore,
$(3)\Rightarrow(2)$ can be concluded.

We next show
 $(2)\Rightarrow (1)$. Applying Lemma \ref{lem2.17}, there exists $T>0$ such that
$\prod_{k=t_{0}}^{t_{0}+T-1}G(k)$ is scrambling for any $t_{0}$.
There exists $\delta>0$ such that $\eta(B(T,t_{0}))>\delta>0$ for
all $t_{0}\ge 0$ because of the compactness of the set ${\mathcal
M}$. So,
\begin{eqnarray}
{\rm diam}(B(t,t_{0}))&=&{\rm diam}\bigg\{B({\rm mod}(t,T),t_{0}+[\frac{t}{T}]T)
\prod\limits_{k=1}^{[\frac{t}{T}]}B(T,t_{0}+(k-1)T)\bigg\}\nonumber\\
&\le& {\rm diam}\bigg\{\prod\limits_{k=1}^{[\frac{t}{T}]}B(t_{0}+k
T-1,t_{0}+(k-1)T)\bigg\}\nonumber\\
 &\le& 2(1-\delta)^{[\frac{t}{T}]}\label{fix}
\end{eqnarray}
holds for any $t_{0}\ge 0$. Here, $[t/T]$ denotes the largest
integer less than $t/{T}$ and ${\rm mod}(t,T)$ denotes the modulus
of the division $t\div T$. Thus,
\begin{eqnarray*}
{\rm diam}({\mathcal G})\le (1-\delta)^{\frac{1}{T}}<1.
\end{eqnarray*}
This proves $(2)\Rightarrow (1)$. Since $(1)\Rightarrow (3)$ is
clear, the theorem is proved.

{\sc Remark 4.} According to Lemma \ref{lem2.17}, it can be seen
that the union of graphs across any time interval of length $T$
has a spanning tree if and only if a union of graphs across any
time interval of length $(m-1)T$ is scrambling.

Moreover, from \cite{Dau}, we conclude more results on the
ergodicity of stochastic matrix sequences as follows:
\begin{proposition}\label{prop4.3}
The implication $(1)\Rightarrow (2)\Rightarrow (3)$ holds for the
following statements:
\begin{enumerate}
\item ${\rm diam}({\mathcal G})<1$;

\item ${\mathcal G}$ is ergodic;

\item for any $t_{0}\ge 0$, the union $\bigcup_{k\ge
t_{0}}\Gamma(k)$ has a spanning tree.
\end{enumerate}
\end{proposition}
{\sc Remark 5.} It should be pointed out that  the implications in
Proposition \ref{prop4.3} cannot be reversed. Counterexamples can
be found in \cite{Mor}. However, in \cite{Mor}, it is also proved
under certain conditions that if the stochastic matrices have the
property that $G_{ij}>0$ if and only if $G_{ji}>0$, then statement
2 is equivalent to statement 3.

Assembling Theorem \ref{thm4.2}, Proposition \ref{prop4.3}, and
the results in \cite{Mor}, it can be shown that, for ${\mathcal
G}\subset{\mathcal M}$, the following implications hold
\begin{eqnarray*}
A_{1}\Leftrightarrow A_{2}\Leftrightarrow A_{3}\Rightarrow
A_{4}\Rightarrow A_{5}
\end{eqnarray*}
regarding the statements:
\begin{itemize}
\item $A_{1}$: ${\rm diam}({\mathcal G})<1$;

\item $A_{2}$: there exists $T>0$ such that the union across any
$T$-length time interval $[t_{0}.t_{0}+T]$:
$\bigcup_{k=t_{0}}^{t_{0}+T}\Gamma(k)$ has a spanning tree;

\item $A_{3}$: ${\mathcal G}$ is uniformly ergodic;

\item $A_{4}$: ${\mathcal G}$ is ergodic;

\item$A_{5}$: for any $t_{0}$, the union across $[t_{0},\infty)$:
$\bigcup_{k\ge t_{0}}\Gamma(k)$ has a spanning tree.
\end{itemize}

In the following, we present some special classes of examples of
 coupled map lattices with time-varying couplings. These
classes were widely used to describe discrete-time networks and
studied in some recent papers \cite{Jost,Ding,Shen,Dau}. The
synchronization criterion for these classes can be verified by
numerical methods. Thus, the synchronizability ${\rm
diam}({\mathcal G})$ of the time-varying couplings can also be
computed numerically.
\subsection{Static topology}
If $G(t)$ is a static matrix, i.e., $G(t)=G$, for all $t\in
{\mathbb Z}^{+}$, then we can write the coupled system
(\ref{lcmls}) as
\begin{eqnarray}
x(t+1)=GF(x(t)).\label{static}
\end{eqnarray}
\begin{proposition}\label{prop4.4}
Let $1=\sigma_{0}$, $\sigma_{1}$, $\sigma_{2}$, $\cdots$,
$\sigma_{m-1}$ be the eigenvalues of $G$ ordered by
$1\ge|\sigma_{1}|\ge|\sigma_{2}|\ge\cdots\ge|\sigma_{m-1}|$. If
$|\sigma_{1}|e^{\mu}<1$, then the coupled system (\ref{static}) is
synchronized.
\end{proposition}
{\em Proof.}  Let $v_{0}=e_{0}$ and choose column vectors $v_{1}$,
$v_{2}$, $\cdots$, $v_{m-1}$ in $\mathbb{R}^{m}$ such that
${v_{0},v_{1},\cdots,v_{m-1}}$ is an orthonormal basis for
$\mathbb{R}^{m}$. Let $A=[v_{0},v_{1},\cdots,v_{m-1}]$. Then,
\begin{eqnarray*}
A^{-1}GA=\left[\begin{array}{ll}1&\alpha\\
0&\hat{G}\end{array}\right],
\end{eqnarray*}
where the eigenvalues of $\hat{G}$ are $\sigma_{1}$, $\cdots$,
$\sigma_{m-1}$. By the Householder theorem (see Theorem 4.2.1 in
\cite{Ser}), for any $\epsilon>0$, there must exist a norm in
$\mathbb{R}^{m}$ such that with its induced matrix norm,
\begin{eqnarray*}
|\sigma_{1}|\le \|\hat{G}\|\le|\sigma_{1}|+\epsilon.
\end{eqnarray*}
Since $\epsilon$ is arbitrary, for the static stochastic matrix
sequence ${\mathcal G}_{0}=\{G,G,\cdots,\}$, it can be concluded
that $\hat{\rho}({\mathcal G}_{0})=|\sigma_{1}|$. Using Theorem
\ref{thm4.1}, the conclusion follows. Moreover, it can be also
obtained that the convergence rate is
$O((|\sigma_{1}|e^{\mu})^{t})$.

{\sc Remark 6.} Similar results have been obtained by several
papers concerning synchronization of coupled map lattices with
static connections (see \cite{Jost,Ding,Lu,Atay-PhysD}). Here, we have proved
this result in a different way as a consequence of our main
result.

\subsection{Finite topology set}
Let ${\mathcal Q}$ be a compact stochastic matrix set
 satisfying {\bf H4}. Consider the following
inclusions:
\begin{eqnarray}
x(t+1)\in {\mathcal Q}F(x(t)),\label{inclusion}
\end{eqnarray}
i.e.,
\begin{eqnarray}
x(t+1)&=&G(t)F(x(t))\label{equation}\\
G(t)&\in&{\mathcal Q}.
\end{eqnarray}
Then the synchronization of the coupled system (\ref{inclusion}) can
be formulated as follows.
\begin{definition}\label{def4.5}
The coupled inclusion system (\ref{inclusion}) is said to be
synchronized if for any stochastic matrix sequence ${\mathcal
G}\subset{\mathcal Q}$, the coupled system (\ref{equation}) is
synchronized.
\end{definition}

In \cite{Shen}, the authors defined the Hajnal diameter and
projection joint spectral radius for a compact stochastic matrix
set.
\begin{definition}\label{def4.6} For the stochastic matrix
set ${\mathcal Q}$, the Hajnal diameter is given by
$${\rm diam}({\mathcal
Q})=\overline{\lim\limits_{t\rightarrow\infty}}\sup\limits_{G(k)\in{\mathcal
Q}}\bigg\{{\rm
diam}(\prod\limits_{k=0}^{t-1}G(k))\bigg\}^{\frac{1}{t}},$$ and
the projection joint spectral radius is
$$\hat{\rho}({\mathcal Q})=\overline{\lim\limits_{t\rightarrow\infty}}
\bigg\{\sup\limits_{G(k)\in{\mathcal Q}}
\|\prod\limits_{k=0}^{t-1}\hat{G}(k)\|\bigg\}^{\frac{1}{t}}.$$
\end{definition}

The following result is from \cite{Shen}.

\begin{lemma} \label{lem4.7}Suppose ${\mathcal Q}$ is a compact set of
stochastic matrices. Then,
\begin{eqnarray*}
{\rm diam}({\mathcal Q})=\hat{\rho}({\mathcal Q}).
\end{eqnarray*}
\end{lemma}

Using Theorem \ref{thm4.1}, we have

\begin{theorem}\label{thm4.8}
If ${\rm diam}({\mathcal Q})e^{\mu}<1$, then the coupled system
(\ref{inclusion}) is synchronized.
\end{theorem}

Moreover, we conclude that the synchronization is uniform with
respect to $t_{0}\in {\mathbb Z}^{+} $ and stochastic matrix
sequences ${\mathcal G}\subset{\mathcal Q}$. Furthermore, we have
the following result on synchronizability of the stochastic matrix
set ${\mathcal Q}$.
\begin{proposition}\label{prop4.9} Let ${\mathcal Q}$ be a compact set of stochastic matrices
satisfying hypothesis {\bf H4}. Then the following statements are
equivalent:
\begin{itemize}

\item ${\mathcal B}_{1}$: ${\rm diam}({\mathcal Q})<1$;

\item ${\mathcal B}_{2}$: for any stochastic matrix sequence
${\mathcal G}\subset{\mathcal Q}$, ${\mathcal G}$ is ergodic;

\item ${\mathcal B}_{3}$: each corresponding graph of a stochastic
matrix $G\in{\mathcal Q}$ has a spanning tree.
\end{itemize}
\end{proposition}

{\em Proof.} The implication ${\mathcal B}_{1}\Rightarrow {\mathcal
B}_{2}\Rightarrow {\mathcal B}_{3}$ is clear by Proposition
\ref{prop4.3}. And ${\mathcal B}_{3}\Rightarrow{\mathcal B}_{1}$ can
be obtained by the proof of Theorem \ref{thm4.2} since ${\mathcal
Q}$ is a finite set of stochastic matrices satisfying hypothesis
{\bf H4}.

{\sc Remark 7.} By the methods introduced in
\cite{Gri,Zhou1,Zhou2}, $\hat{\rho}({\mathcal Q})$ can be computed
to arbitrary precision for a finite set ${\mathcal Q}$ despite a
large computational complexity.

\subsection{Multiplicative ergodic topology sequence}
Consider the stochastic matrix sequence ${\mathcal
G}=\{G(t)\}_{t\in{\mathbb Z}^{+}}$ driven by some dynamical system
$\mathcal Y=\{\Omega,{\mathcal F},P,\theta^{(t)}\}$, i.e.,
${\mathcal G}=\{G(\theta^{(t)}\omega)\}$ for some continuous map
$G(\cdot)$. Recall the Lyapunov exponent for ${\mathcal G}$:
\begin{eqnarray*}
\sigma(v,\omega)=\overline{\lim\limits_{t\rightarrow\infty}}\frac{1}{t}
\log\bigg\|\prod\limits_{k=0}^{t-1}G(\theta^{(k)}\omega)v\bigg\|.
\end{eqnarray*}
It is clear that $\sigma(e_{0},\omega)=0$ for all $\omega$ and
$\sigma(v,\omega)\le 0$
 for all $\omega$ and $v\in
\mathbb{R}^{m}$. So, the linear subspace
\begin{eqnarray*}
L_{\omega}=\{v,\sigma(v,\omega)<0\}
\end{eqnarray*}
denotes the directions transverse to the synchronization manifold.
If $P$ is an ergodic measure for the MDS $\mathcal Y$, then
$\sigma(u,\omega)$ and $L(\omega)$ are the same for almost all
$\omega$ with respect to $P$ \cite{Arn}. Then we can let
$\sigma_{1}$ be the largest Lyapunov exponent of ${\mathcal G}$
transverse to the synchronization direction $e_{0}$. By Theorem
\ref{thm4.1} and Corollary \ref{col3.3}, we have
\begin{theorem}\label{thm4.10}
Suppose that $\theta^{(t)}$ is a continuous semiflow, $G(\cdot)$
is continuous on all $\omega\in\Omega$ and non-singular, and
$\Omega$ is compact. If
$$\sup\limits_{\rm{Erf}_{\theta}(\Omega)}\sigma_{1}+\mu<0,$$
 then the coupled system (\ref{lcmls}) is synchronized.
\end{theorem}

{\sc Remark 8.} There are many papers discussing the computation
of multiplicative Lyapunov exponents; for example, see \cite{Vle}.
In particular, \cite{Mai} discussed the Lyapunov exponents for the
product of infinite matrices. By Lemma \ref{lem2.7}, we can
compute the largest projection Lyapunov exponent which equals
$\sigma_{1}$. We will illustrate this in the following section.

\section{Numerical illustrations}
In this section, we will numerically illustrate the theoretical
results on synchronization of CML with time-varying couplings. In
these examples, the coupling matrices are driven by random
dynamical systems which can be regarded as stochastic processes.
Then the projection Lyapunov exponents are be computed numerically
by the time series of coupling matrices. In this way, we can
verify the synchronization criterion and analyze synchronizability
numerically. Consider the following coupled map network with
time-varying topology:
\begin{eqnarray}
x^{i}(t+1)=\frac{1}{\sum\limits_{k=1}^{m}A_{ik}(t)}\sum\limits_{j=1}^{m}A_{ij}(t)f(x^{j}(t)),~i=\onetom,\label{eg}
\end{eqnarray}
where $x^{i}(t)\in \mathbb{R}$ and $f(s)=\alpha s(1-s)$ is the
logistic map with $\alpha=3.9$, which implies that the Lyapunov
exponent of $f$ is $\mu\approx 0.5$. The stochastic coupling
matrix at time $t$ is
$$G(t)=[G_{ij}(t)]_{i,j=1}^{m}=\bigg[\frac{A_{ij}(t)}{\sum\limits_{j=1}^{m}A_{ij}(t)}\bigg]_{i,j=1}^{m}.$$

\subsection{Blinking scale-free networks}
The blinking scale-free network is a model initiated by a
scale-free network and evolves with malfunction and recovery. At
time $t=0$, the initial graph $\Gamma(0)$ is a scale-free network
introduced in \cite{Bar}. At each time $t\ge 1$, every vertex $i$
malfunctions with probability $p\ll 1$. If vertex $i$
malfunctions, all edges linked to it disappear. In addition, a
malfunctioned vertex recovers after a time interval $T$ and then
causes the re-establishment of all edges linked to it in the
initial graph $\Gamma(0)$. The coupling $A_{ij}(t)=A_{ji}(t)=1$ if
vertex $j$ is connected to $i$ at time $t$; otherwise,
$A_{ij}(t)=A_{ji}(t)=0$, and $A_{ii}(t)=1$, for all $i,j=\onetom$.
 \begin{figure}
 \epsfig{file=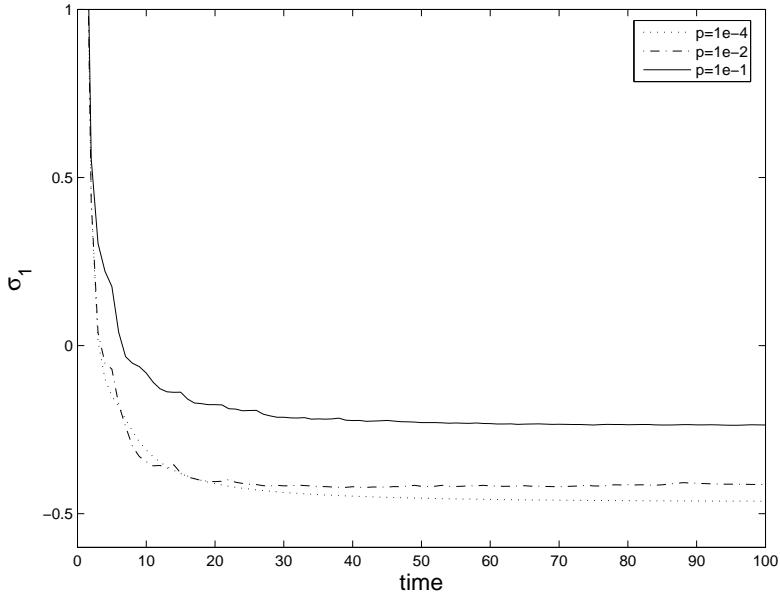,height=8cm} \caption{Convergence of the
 second Lyapunov exponent $\sigma_{1}$ for the blinking topology during the topology
 evolution with the same recovery time $T=3$ and different
 malfunction probability $p=10^{-1}$, $p=10^{-2}$, and $p=10^{-4}$.
 The initial scale-free graph is constructed by the method
 introduced in \cite{Bar} with network size $500$ and average
 degree $12$.}\label{blinking1}
 \end{figure}

\begin{figure}
 \epsfig{file=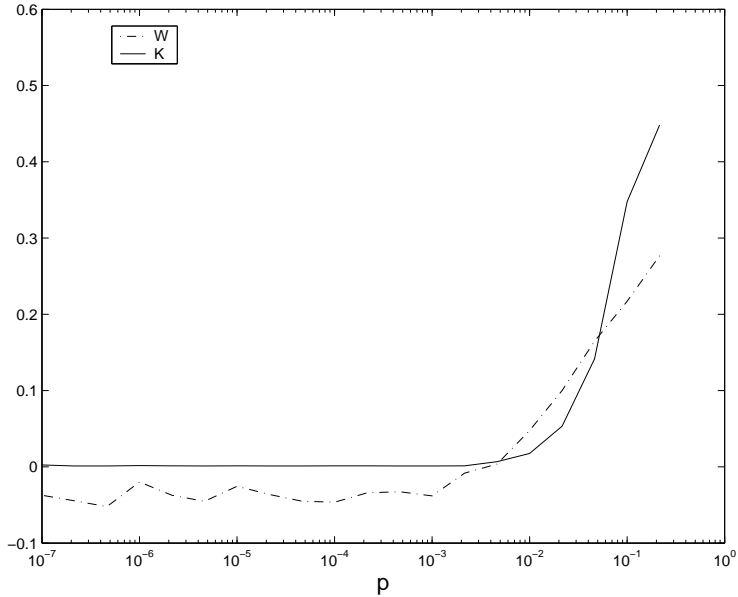,height=8cm} \caption{Variation of $K$
 and $W$ with respect to $p$ for the blinking topology.}\label{blinking2}
 \end{figure}

In Figure \ref{blinking1}, we show the convergence of the second
Lyapunov exponent $\sigma_{1}$ during the topology evolution with
different malfunction probability $p$. We measure synchronization
by the variance
$K=1/(m-1)<\sum_{i=1}^{m}(x^{i}(t)-\bar{x}(t))^{2}>$, where
$<\cdot>$ denotes the time average, and denote $W=\sigma_{1}+\mu$.
We pick the evolution time length to be $1000$ and choose initial
conditions randomly from the interval $(0,1)$. In Figure
\ref{blinking2}, we show the variation of $K$ and $W$ with respect
to the malfunction probability $p$. It can be seen that the region
where $W$ is negative coincides with the region of
synchronization, i.e., where $K$ is near zero.

\subsection{Blurring directed graph process}
A blurring directed graph process is one where each edge weight is
a modified Wiener process. In details, the graph process is
started with a directed weighted graph $\Gamma(0)$ of which for
each vertex pair $(i,j)$, one of two edges $A_{ij}(0)$ and
$A_{ji}(0)$ is a random variable uniformly distributed between $1$
and $2$, and the other is zero with equal probability, for all
$i\ne j$; $A_{ii}(0)=0$ for all $i=\onetom$. At each time $t\ge
1$, for each $A_{ij}(t-1)\ne 0$, $i\ne j$ we denote the difference
$A_{ij}(t)-A_{ij}(t-1)$ by a Gaussian distribution ${\mathcal
N}(0,r^2)$ which is statistically independent for all $i\ne j$ and
$t\in\mathbb Z^{+}$. If resulted in that $A_{ij}(t)$ is negative,
a weight will be added to the reversal orientation, i.e.,
$A_{ji}(t)=|A_{ij}(t)|$ and $A_{ij}(t)=0$. Moreover, if the
process above results in that there exists some index $i$ such
that $A_{ij}=0$ holds for all $j=\onetom$, then pick
$A_{ii}(t)=1$.

 \begin{figure}
 \epsfig{file=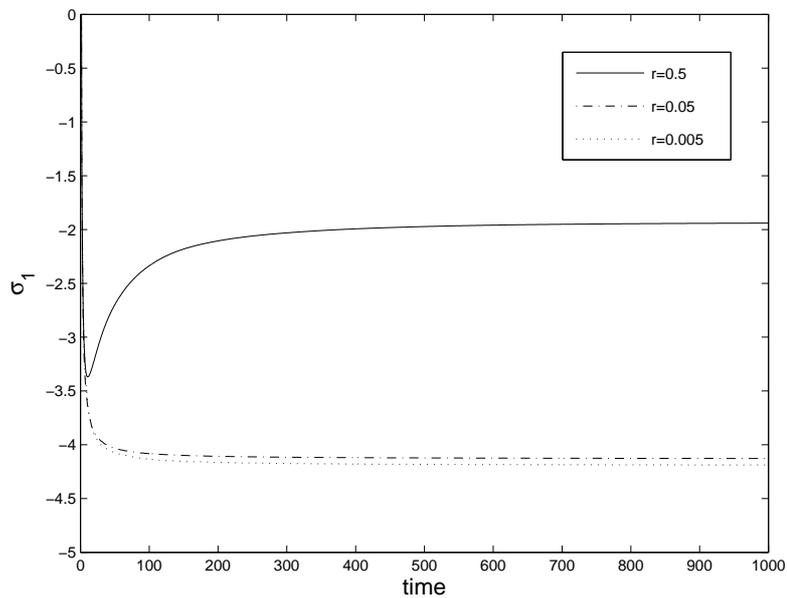,height=8cm} \caption{Convergence of the
 second Lyapunov exponent $\sigma_{1}$ for the blurring graph process during the topology
 evolution with Gaussian variance $r=0.5$, $0.05$,
 $0.005$, and the size of the network $m=100$.}\label{blur1}
 \end{figure}

\begin{figure}
 \epsfig{file=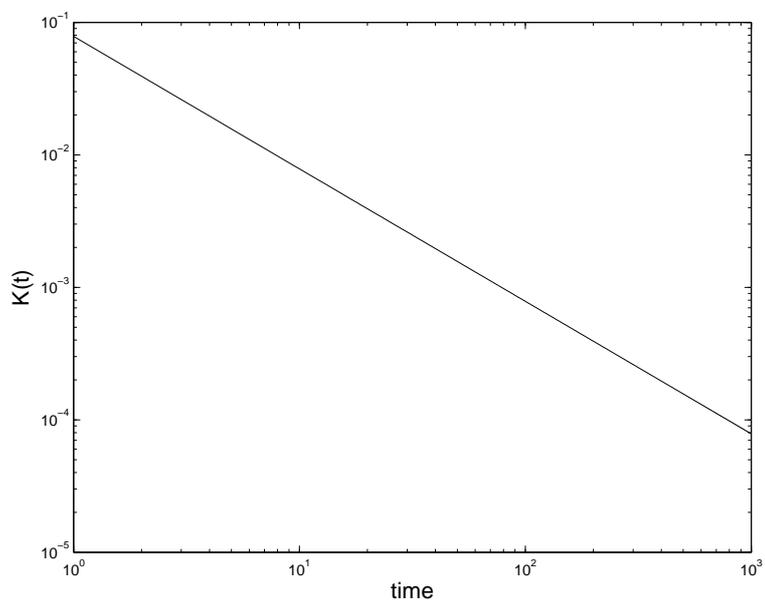,height=8cm} \caption{Variation of $K(t)$
 with respect to time for the blurring graph process.}\label{blur2}
 \end{figure}

In Figure \ref{blur1}, we show the convergence of the second
Lyapunov exponent $\sigma_{1}$ during the topology evolution for
different values of the Gaussian distribution variance $r$.
Picking $r=0.05$, we show the synchronization of the coupled
system (\ref{eg}). Let
$K(t)=1/(m-1)<\sum_{i=1}^{m}(x^{i}(t)-\bar{x}(t))^{2}>_{t}$, where
$<\cdot>_{t}$ denotes the time average from $0$ to $t$. Since
$W=\sigma_{1}+\mu$ is about $-0.6$, i.e. less than zero, the
coupled system is synchronized. Figure \ref{blur2} shows in
logarithmic scale the convergence of $K(t)$ to zero.

\section{Conclusion}
In this paper, we have presented a synchronization analysis for
discrete-time dynamical networks with time-varying topologies. We
have extended the concept of the Hajnal diameter to generalized
matrix sequences to discuss the synchronization of the coupled
system. Furthermore, this quantity is equivalent to other widely
used quantities such as the projection joint spectral radius and
transverse Lyapunov exponents, which we have also extended to the
time-varying case. Thus, these results can be used to discuss the
synchronization of the CML with time-varying couplings. The Hajnal
diameter is utilized to describe synchronizability of the
time-varying couplings and obtain a criterion guaranteeing
synchronization. Time-varying couplings can be regarded as a
stochastic matrix sequence associated with a sequence of graphs.
Synchronizability is tightly related to the topology. As we have
shown, the statement that ${\rm diam}({\mathcal G})<1$, i.e. that
chaotic synchronization is possible, is equivalent to saying that
there exists an integer $T$ such that the union of the graphs
across any time interval of length $T$ has a spanning tree. The
methodology will be similarly extended to higher dimensional maps
elsewhere.

\section*{Appendix}
{\em Proof.} ({\bf Lemma \ref{lem2.4}}) The proof of this lemma
comes from \cite{Shen} with a minor modification. First, we show
${\rm diam}({\mathcal L},\phi)\le\hat{\rho}({\mathcal L},\phi)$. Let
$J$ be any complement of $\mathcal E_{0}$ in $\mathbb{R}^{m}$ with a
basis $u_{0},\cdots,u_{m-1}$ such that $u_{0}=e_{0}$. Let
$A=[u_{0},u_{1},\cdots,u_{m-1}]$ which is nonsingular. Then, for any
$t>t_{0}$ and $t_{0}\ge 0$,
\begin{eqnarray*}
A^{-1}L_{t}(\varrho^{(t-t_{0})}\phi)A=\left[\begin{array}{cc}c(t)&\alpha_{t}\\
0&\hat{L}_{t}(\varrho^{(t-t_{0})}\phi)\end{array}\right],
\end{eqnarray*}
where $c(t)$ denotes the row sum of
$L_{t}(\varrho^{(t-t_{0})}\phi)$ which is also the eigenvalue
corresponding eigenvector $e$ and
$\hat{L}_{t}(\varrho^{(t-t_{0})}\phi)$ can be the solution of
linear equation (\ref{proj}) with $P$ composed of the rows of
$A^{-1}$ except the first row. For any $d>\hat{\rho}({\mathcal
L},\phi)$, there exists $T>0$ such that the inequality
\begin{eqnarray*}
\bigg\|\prod\limits_{k=t_{0}}^{t_{0}+t-1}\hat{L}_{k}(\varrho^{(k-t_{0})}\phi)\bigg\|\le
d^{t}
\end{eqnarray*}
holds for all $t\ge T$ and $t_{0}\ge 0$. Let
\begin{eqnarray*}
A^{-1}\prod\limits_{k=t_{0}}^{t_{0}+t-1}L_{k}(\varrho^{(k-t_{0})}\phi)A=\left[\begin{array}{cc}
\prod\limits_{k=t_{0}}^{t_{0}+t-1}c(k)&\alpha_{t}\\
0&\prod\limits_{k=t_{0}}^{t_{0}+t-1}\hat{L}_{k}(\varrho^{(k-t_{0})}\phi)\end{array}\right]
\end{eqnarray*}
Then,
\begin{eqnarray*}
\bigg\|A^{-1}\prod\limits_{k=t_{0}}^{t_{0}+t-1}L_{k}(\varrho^{(k-t_{0})}\phi)A-
\left[\begin{array}{c}1\\0\\\vdots\\0\end{array}\right]
\big(\prod\limits_{k=t_{0}}^{t_{0}+t-1}c(k),\alpha_{t}\big)\bigg\|=\bigg\|\left[\begin{array}{cc}0&0\\
0&\prod\limits_{k=t_{0}}^{t_{0}+t-1}\hat{L}_{k}(\varrho^{(k-t_{0})}\phi)\end{array}\right]\bigg\|\le
Cd^{t}
\end{eqnarray*}
holds for some constant $C>0$. Therefore,
\begin{eqnarray*}
\bigg\|\prod\limits_{s=t_{0}}^{t_{0}+t-1}L_{k}(\varrho^{(k-t_{0})}\phi)-A
\left[\begin{array}{c}1\\0\\\vdots\\0\end{array}\right]
\big(\prod\limits_{k=t_{0}}^{t_{0}+t-1}c(k),\alpha_{t}\big)A^{-1}\bigg\|&\le&
C_{1}d^{t},\\
\bigg\|\prod\limits_{k=t_{0}}^{t_{0}+t-1}L_{k}(\varrho^{(k-t_{0})}\phi)-e\cdot
q\bigg\|&\le& C_{1}d^{t},
\end{eqnarray*}
where
$q=\big[\prod_{k=t_{0}}^{t_{0}+t-1}c(k),\alpha_{t}\big]A^{-1}$ and
$C_{1}$ is a positive constant. It says that all row vectors of
$\prod_{k=t_{0}}^{t_{0}+t-1}L_{k}(\varrho^{(k-t_{0})}\phi)$ lie
inside the $C_{1}d^{m}$ neighborhood of $q$. Hence,
\begin{eqnarray*}
{\rm
diam}\bigg(\prod\limits_{k=t_{0}}^{t_{0}+t-1}L_{k}(\varrho^{(k-t_{0})}\phi)\bigg)\le
C_{2}d^{t}
\end{eqnarray*}
for some constant $C_{2}>0$, all $t\ge T$, and $t_{0}\ge 0$. This
implies that ${\rm diam}({\mathcal L},\phi)\le d$. Since $d$ is
arbitrary, ${\rm diam}({\mathcal L},\phi)\le\hat{\rho}({\mathcal
L},\phi)$ can be concluded.

Second, we show that $\hat{\rho}({\mathcal L},\phi)\le {\rm
diam}({\mathcal L},\phi)$. For any $d>{\rm diam}({\mathcal
L},\phi)$, there exists $T>0$ such that
\begin{eqnarray*}
{\rm
diam}\bigg(\prod\limits_{k=t_{0}}^{t_{0}+t-1}L_{k}(\varrho^{(k-t_{0})}\phi)\bigg)\le
d^{t}.
\end{eqnarray*}
holds for all $t\ge T$ and $t_{0}\ge 0$. Letting $q$ be the first
row of
$\prod_{k=t_{0}}^{t_{0}+t-1}L_{k}(\varrho^{(k-t_{0})}\phi)$, we
have
\begin{eqnarray*}
\left\|\prod\limits_{k=t_{0}}^{t_{0}+t-1}L_{k}(\varrho^{(k-t_{0})}\phi)-e\cdot
q\right\|\le C_{3}d^{t}
\end{eqnarray*}
for some positive constant $C_{3}$. Let $A$ be defined as above.
Then,
\begin{eqnarray*}
\left\|A^{-1}\prod\limits_{k=t_{0}}^{t_{0}+t-1}L_{k}(\varrho^{(k-t_{0})}\phi)A-A^{-1}e\cdot
q A\right\|\le C_{4} d^{t},
\end{eqnarray*}
i.e.,
\begin{eqnarray*}
\left\|\left[\begin{array}{cc}\prod\limits_{k=t_{0}}^{t_{0}+t-1}c(k)&\alpha_{t}\\
0&\prod\limits_{k=t_{0}}^{t_{0}+t-1}\hat{L}_{k}(\varrho^{(k-t_{0})}\phi)
\end{array}\right]-\left[\begin{array}{cc}\gamma&\beta\\
0&0\end{array}\right]\right\|\le C_{4}d^{t}
\end{eqnarray*}
holds for some $\gamma$ and $\beta$. This implies that
\begin{eqnarray*}
\left\|\prod\limits_{k=t_{0}}^{t_{0}+t-1}\hat{L}_{k}(\varrho^{(k-t_{0})}\phi)\right\|\le
C_{5}d^{t}
\end{eqnarray*}
holds for all $t\ge T$, $t_{0}\ge 0$, and some $C_{5}>0$. Therefore,
$\hat{\rho}({\mathcal L},\phi)\le d$. The proof is completed since
$d$ is chosen arbitrarily.

{\em Proof.} ({\bf Lemma \ref{lem2.6}}) Let
$\hat{\lambda}_{\max}=\sup_{v\in
\mathbb{R}^{m-1}}\hat{\lambda}({\mathcal L},\phi,v)$. First, it is
easy to see that $\log\hat{\rho}({\mathcal L},\phi)\ge
\hat{\lambda}_{\max}$. We will show $\log\hat{\rho}({\mathcal
L},\phi)=\hat{\lambda}_{\max}$. Otherwise, there exists
$d\in(\exp(\hat{\lambda}_{\max}),\hat{\rho}({\mathcal L},\phi))$. By
the properties of Lyapunov exponents, for any normalized orthogonal
basis $u_{1},u_{2},\cdots,u_{m-1}\in \mathbb{R}^{m-1}$ with Lyapunov
exponent $\hat{\lambda}({\mathcal L},\phi,u_{i})=\hat{\lambda}_{i}$,
then for any $u\in \mathbb{R}^{m-1}$ we have
$\hat{\lambda}({\mathcal L},\phi,u)=\hat{\lambda}_{i_{u}}$, where
$i_{u}\in\{1,2,\cdots,m-1\}$. $\hat{\rho}({\mathcal L},\phi)>d$
implies that there exist $t_{0}\ge 0$ and a sequence $t_{n}$ with
$\lim_{n\rightarrow\infty}t_{n}=+\infty$ such that
\begin{eqnarray*}
\bigg\|\prod\limits_{k=t_{0}}^{t_{n}+t_{0}-1}\hat{L}_{k}(\varrho^{(k-t_{0})}\phi)\bigg\|>d^{t_{n}}
\end{eqnarray*}
for all $n\ge 0$. That is, there also exists a sequence $v_{n}\in
\mathbb{R}^{m-1}$ with $\|v_{n}\|=1$ such that
\begin{eqnarray*}
\bigg\|\prod\limits_{k=t_{0}}^{t_{n}+t_{0}-1}\hat{L}_{k}(\varrho^{(k-t_{0})}\phi)v_{n}\bigg\|>d^{t_{n}}.
\end{eqnarray*}
There exists a subsequence of $v_{n}$ (still denoted by $v_{n}$)
with $\lim_{n\rightarrow\infty}v_{n}=v^{*}$. Let $\delta
v_{n}=v_{n}-v^{*}$. We have
\begin{eqnarray*}
\bigg\|\prod\limits_{k=t_{0}}^{t_{n}+t_{0}-1}\hat{L}_{k}(\varrho^{(k-t_{0})}\phi)v^{*}\bigg\|\ge
\bigg\|\prod\limits_{k=t_{0}}^{t_{n}+t_{0}-1}\hat{L}_{k}(\varrho^{(k-t_{0})}\phi)v_{n}\bigg\|
-\bigg\|\prod\limits_{k=t_{0}}^{t_{n}+t_{0}-1}\hat{L}_{k}(\varrho^{(k-t_{0})}\phi)\delta
v_{n}\bigg\|.
\end{eqnarray*}
Note that we can write $\delta v_{n}=\sum_{i=1}^{m-1}\delta
x^{i}_{n} u_{i}$ where $\delta x^{i}_{n}\in \mathbb{R}$ with
$\lim_{n\rightarrow\infty}\delta x^{i}_{n}=0$. So, there exists an
integer $N$ such that
$\big\|\prod_{k=t_{0}}^{t_{n}+t_{0}-1}\hat{L}_{k}(\varrho^{(k-t_{0})}\phi)\delta
v_{n}\big\|\le \big(\sum_{i=1}^{m-1}|\delta
x_{n}^{i}|\big)d^{t_{n}}$ holds for all $n\ge N$. Then, we have
\begin{eqnarray*}
\bigg\|\prod\limits_{k=t_{0}}^{t_{n}+t_{0}-1}\hat{L}_{k}(\varrho^{(k-t_{0})}\phi)v^{*}\bigg\|\ge
d^{t_{n}} -d^{t_{n}}\bigg(\sum\limits_{i=1}^{m-1}|\delta
x_{n}^{i}|\bigg)\ge C d^{t_{n}} .
\end{eqnarray*}
for all $n\ge N$ and some $C>0$. This implies $\max_{v\in
\mathbb{R}^{m}}\hat{\lambda}({\mathcal L},\phi,v)\ge \log d$ which
contradicts with the assumption
$d\in(\exp(\hat{\lambda}_{\max}),\hat{\rho}(\mathcal L,\phi))$.
Hence, $\hat{\lambda}_{\max}=\log\hat{\rho}({\mathcal L},\phi)$.

{\em Proof.} ({\bf Lemma \ref{lem2.7}}) Recalling that
$\{\Phi,\mathcal B,P,\varrho^{(t)}\}$ denotes a random dynamical
system, where $\Phi$ denotes the state space, $\mathcal B$ denotes
the $\sigma$-algebra, $P$ denotes the probability measure, and
$\varrho^{(t)}$ denotes the semiflow. For a given $\phi\in\Phi$ we
denote $L(\varrho^{(t)}\phi)$ by $L(t)$. Let
$A=[u_{1},u_{2},\cdots,u_{m}]\in \mathbb{R}^{m\times m}$ where
$u_{1},\cdots,u_{m}$ denotes a basis of $\mathbb{R}^{m}$ and
$u_{1}=e$,
$$A^{-1}=\left[\begin{array}{c}v_{1}\\v_{2}\\\vdots\\v_{m}\end{array}\right]\in
\mathbb{R}^{m\times m}$$ is the inverse of $A$ with
\begin{eqnarray*}
\bar{L}(t)=A^{-1}L(t)A=\left[\begin{array}{cc}c(t)&\alpha^{\top}(t)\\
0&\hat{L}(t)\end{array}\right],~\hat{L}(t)=A^{*}_{1}D(t)A_{1},~\alpha^{\top}(t)=v_{1}L(t)A_{1},
\end{eqnarray*}
where $A_{1}=[u_{2},\cdots,u_{m}]\in \mathbb{R}^{m\times (m-1)}$ and
$$A^{*}_{1}=\left[\begin{array}{c}v_{2}\\\vdots\\v_{m}\end{array}\right]\in
\mathbb{R}^{(m-1)\times m}.$$ One can see that the set of Lyapunov
exponents of the dynamical system $\{\bar{L}(t)\}_{t\in\mathbb
Z^{+}}$ are the same as those of $\{L(t)\}_{t\in\mathbb Z^{+}}$.
For any $z(0)=\left[x(0),y(0)\right]\in \mathbb{R}^{m}$ where
$x(0)\in \mathbb{R}$ and $y(0)\in \mathbb{R}^{m-1}$, this
evolution $z(t+1)=\bar{L}(t)z(t)$ leads
\begin{eqnarray*}
z(t)=\left[\begin{array}{l}x(t)\\y(t)\end{array}\right]=\left[\begin{array}{l}
c(t-1)x(t-1)
+\alpha^{\top}(t-1)y(t-1)\\\hat{L}(t-1)y(t-1)\end{array}\right].
\end{eqnarray*}
So, we have
\begin{eqnarray}
y(t)&=&\prod\limits_{k=0}^{t-1}\hat{L}(k)y(0)\nonumber\\
x(t)&=&\prod\limits_{k=0}^{t-1}c(k)x(0)+\sum\limits_{k=1}^{t}
\prod\limits_{p=t-k+1}^{t-1}c(p)
\alpha^{\top}(t-k)\prod\limits_{q=0}^{t-k-1}\hat{L}(q)y(0)\nonumber\\
&&\label{xy}
\end{eqnarray}
If the upper bound is less than the lower bound for the left
matrix product $\prod$, then the product should be the identity
matrix. In the following, we denote by $\hat{\mathcal L}$ the
projection sequence map of ${\mathcal L}$ and will prove this
lemma for two cases.

{\bf\em Case 1}: $\lambda_{0}\le\log\hat{\rho}({\mathcal
L},\phi)$. Since $\hat{\rho}({\mathcal L},\phi)$ is just the
largest Lyapunov exponent of $\hat{\mathcal L}$ defined by
$\hat{\lambda}$, from conditions 1 and 2, one can see that for any
$\epsilon>0$, there exists $T>0$ such that for any $t\ge T$, it
holds that $|\alpha(t)|\le e^{\epsilon t}$,
$\big\|\prod_{k=0}^{t-1}\hat{L}(k)\big\|\le
e^{(\hat{\lambda}+\epsilon)t}$, and
$e^{(\lambda_{0}-\epsilon)t}\le|\prod_{k=0}^{t-1}c(k)|\le
e^{(\lambda_{0}+\epsilon)t}$. Thus, we can obtain
\begin{eqnarray*}
&&\prod\limits_{k=t-k+1}^{t-1}|c(p)|=\prod\limits_{p=0}^{t-1}|c(p)|\times\frac{1}{\prod\limits_{p=0}^{t-k}|c(p)|}\\
&&=\left\{\begin{array}{ll}e^{(\lambda_{0}+\epsilon)(t)}e^{-(\lambda_{0}-\epsilon)(t-k+1)}&k\le
t-T+1,\\
e^{(\lambda_{0}+\epsilon)(t-1)}\max\limits_{T\ge q\ge
0}\bigg(\prod\limits_{p=0}^{q}|c(p)|\bigg)^{-1}&t-1\ge k\ge
t-T.\end{array}\right.
\end{eqnarray*}
Then, we have
\begin{eqnarray*}
|x(t)|&\le&
\prod\limits_{k=0}^{t-1}|c(k)||x(0)|+\sum\limits_{k=1}^{t-T+1}
\prod\limits_{p=t-k+1}^{t-1}|c(p)|
|\alpha^{\top}(t-k)|\prod\limits_{q=0}^{t-k-1}\|\hat{L}(q)\|\|y(0)\|\\
&&+\sum\limits_{k=t-T}^{t-1}\prod\limits_{p=t-k+1}^{t-1}|c(p)|
\|\alpha^{\top}(t-k)\|\prod\limits_{q=0}^{t-k-1}\|\hat{L}(q)\|\|y(0)\|\\
&\le&
e^{(\lambda_{0}+\epsilon)t}+\sum\limits_{k=1}^{t-T+1}e^{(\lambda_{0}+\epsilon)(t-1)}e^{\epsilon
t} e^{-(\lambda_{0}-\epsilon)(t-k)}
e^{(\hat{\lambda}+\epsilon)(t-k)}
+M_{1}e^{(\lambda_{0}+\epsilon)(t-1)}\\
&\le&e^{(\hat{\lambda}+\epsilon)t}+e^{(\hat{\lambda}+4\epsilon)t}e^{-(\lambda_{0}+\epsilon)}\sum\limits_{k=1}^{t-T+1}
e^{(-\hat{\lambda}+\lambda_{0}-3\epsilon)k}+M_{1}e^{(\lambda_{0}+\epsilon)t}\\
&\le&M_{2}e^{(\hat{\lambda}+4\epsilon)t},
\end{eqnarray*}
where $$M_{1}=(T+1)\max\limits_{T\ge q\ge
0}\bigg(\prod\limits_{p=0}^{q}|c(p)|\bigg)^{-1}e^{\epsilon
T}\bigg(\prod\limits_{p=0}^{q}\|\hat{L}(p)\|\bigg)\|y(0)\|$$ $$
M_{2}=1+M_{1}+e^{-(\lambda_{0}+\epsilon)}\sum\limits_{k=1}^{\infty}
e^{-3\epsilon k}.$$ So,
\begin{eqnarray*}
\overline{\lim_{t\rightarrow\infty}}\frac{1}{t}\log\|z(t-1)\|\le\hat{\lambda}+4\epsilon
\end{eqnarray*}
holds for all $z(0)\in \mathbb{R}^{m}$. Noting that
$\hat{\lambda}$ must be less than the largest Lyapunov exponent of
${\mathcal L}$, we conclude that $\hat{\lambda}$ is right the
largest Lyapunov exponent. This implies the conclusion of the
lemma.

{\bf\em Case 2}: $\lambda_{0}>\hat{\lambda}$. Noting that for any
$\epsilon\in(0,(\lambda_{0}-\hat{\lambda})/3)$, there exists $T$
such that
\begin{eqnarray}
\prod\limits_{k=0}^{t}|c^{-1}(k)|\|\alpha^{\top}(t)\|\prod\limits_{l=0}^{t}
\|\hat{L}(l)\|\le C
e^{(-\lambda_{0}+\hat{\lambda}+3\epsilon)t}\label{lim}
\end{eqnarray}
for all $t\ge T$ and some constant $C>0$. Let
\begin{eqnarray*}
x=-\sum\limits_{t=0}^{\infty}\prod\limits_{k=0}^{t}c^{-1}(k)
\alpha^{\top}(t)\prod\limits_{l=0}^{t-1} \hat{L}(l) y,
\end{eqnarray*}
which in fact exists and is finite according to the inequality
(\ref{lim}). Then, let
\begin{eqnarray*}
V_{\phi}=\bigg\{z=\left[\begin{array}{c}x\\y\end{array}\right]:~
x+\sum\limits_{t=0}^{\infty}\prod\limits_{k=0}^{t}c^{-1}(k)
\alpha^{\top}(t)\prod\limits_{l=0}^{t-1} \hat{L}(l) y=0\bigg\}
\end{eqnarray*}
be the transverse space. For any
$\left[\begin{array}{c}x(0)\\y(0)\end{array}\right]\in V_{\phi}$,
\begin{eqnarray*}
x(t)=-\sum\limits_{k=t}^{\infty}\prod\limits_{p=t}^{k}
c^{-1}(p)\alpha^{\top}(k)\prod\limits_{q=0}^{k-1} \hat{L}(q) y(0).
\end{eqnarray*}
Noting that there exists $T>0$ such that $
\prod\limits_{p=t}^{k}|c^{-1}(p)|\le
e^{(-\lambda_{0}+\epsilon)(k-t)+2\epsilon t}$ for all $t\ge T$, we
have
\begin{eqnarray*}
|x(t)|&\le&\sum\limits_{k=t}^{\infty}\prod\limits_{p=t}^{k}|c^{-1}(p)|\|\alpha^{\top}(k)\|
\bigg\|\prod\limits_{q=0}^{k-1}\hat{L}(q)\bigg\|\|
y(0)\|\\
&\le&\sum\limits_{k=t}^{\infty}e^{(-\lambda_{0}+\epsilon)(k-t)}e^{2\epsilon
t}e^{\epsilon k}e^{(\hat{\lambda}+\epsilon)k}\\
&\le&
\bigg\{\sum\limits_{k=t}^{\infty}e^{(-\lambda_{0}+\hat{\lambda}+3\epsilon)(k-t)}\bigg\}
e^{(\hat{\lambda}+4\epsilon)t}\le M_{2}
e^{(\hat{\lambda}+4\epsilon)t}
\end{eqnarray*}
for all $t\ge T$and some constants $M_{2}>0$. So, it can be
concluded that
\begin{eqnarray*}
\overline{\lim_{t\rightarrow\infty}}\frac{1}{t}\log\|z(t-1)\|\le
\hat{\lambda}+4\epsilon.
\end{eqnarray*}
Since $\epsilon$ is chosen arbitrarily, there exists an $m-1$
dimensional subspace $V_{\phi}=\{z=[x\,\, y]^\top
:~x=-\sum_{t=0}^{\infty}\prod_{k=0}^{t}c^{-1}(k)\alpha^{\top}(t)
\prod_{l=0}^{t-1} \hat{L}(l) y\}$ of which the largest Lyapunov
exponent is less than $\hat{\lambda}$.  The largest Lyapunov
exponent of $V_{\phi}$ is clearly greater than $\hat{\lambda}$.
Therefore, we conclude that $\hat{\lambda}$ i.e.
$\log(\hat{\rho}(L))$, is the largest Lyapunov exponent of $L$
except $\lambda_{0}$. The proof is completed.


\begin{thebibliography}{99}

\bibitem{Ash}
{\sc P. Ashwin, J. Buescu, I. Stewart}, {\em From attractor to
chaotic saddle: a tale of transverse instability}, Nonlinearity, 9
(1996), pp. 703--737.


\bibitem{Mil1}
{\sc J. Milnor}, {\em On the concept of attractors}, Commun. Math.
Phys., 99 (1985), pp. 177--195.


\bibitem{Pik}
{\sc A. Pikovsky, M. Rosenblum, J. Kurths}, {\em Synchronization:
A universal concept in nonlinear sciences}, Cambridge University
Press, 2001.

\bibitem{Kaneko}
{\sc K. Kaneko}, {\em Theory and applications of coupled map lattices},
Wiley, 1993.

\bibitem{Jost}
{\sc J. Jost, M. P. Joy}, {\em Spectral properties and
synchronization in coupled map lattices},  Phys. Rev. E, 65
(2001), 016201.

\bibitem{Ding}
{\sc Y. H. Chen, G. Rangarajan, M. Ding}, {\em General stability
analysis of synchronized dynamics in coupled systems}, Phys. Rev.
E, 67 (2003), 026209.

\bibitem{Atay}
{\sc F. M. Atay, J. Jost, A. Wende}, {\em Delays, connection
topology, and synchronization of coupled chaotic maps}, Phys. Rev.
Lett., 92:14 (2004), 144101.

\bibitem{Lu}
{\sc W. Lu, T. Chen}, {\em Synchronization of linearly coupled
networks with discrete time systems}, Physica D, 198 (2004), pp.
148--168.


\bibitem{Wu}
{\sc C. W. Wu}, {\em Synchronization in networks of nonlinear
dynamical systems coupled via a directed graph}, Nonlinearity, 18
(2005), pp. 1057--1064.

\bibitem{Lu1} {\sc W. Lu, T Chen}, {\em Global synchronization of linearly
coupled Lipschitz map lattices with a directed
graph}, IEEE Transactions on Circuits and Systems-II: Express
Briefs, 54:2 (2007), pp. 136--140.

\bibitem{Sab}
{\sc R. Olfati-Saber, R. M. Murray}, {\em Consensus problems in
networks of agents with switching topology and time-delys}, IEEE
T. Autom. Cont., 49:9 (2004), pp. 1520--1533.



\bibitem{Hat}
{\sc Y.~Hatano, M.~Mesbahi}, {\em Agreement over random networks},
IEEE Conference on Decision and Control. (2004)
http://www.aa.washington.edu/faculty/mesbahi/papers/random-cdc.pdf.


\bibitem{Vic}
{\sc T. Vicsek, A. Czir$\acute{o}$k, E. Ben-Jacob, I. Cohen, O.
Schochet}, {\em Novel type of phase transitions in a system of
self-driven particles}, Phys. Rev. Lett, 75:6 (1995), pp.
1226--1229.


\bibitem{Mor}
{\sc L. Moreau}, {\em Stability of multiagent systems with time
dependent communication links}, IEEE Trans. Auto. Cont., 50:2
(2005), pp. 169--182.

\bibitem{Lv}
{\sc J. H. Lu, G. Chen}, {\em A time-varying complex dynamical
network model and its controlled synchronization criterion}, IEEE
Trans. Auto. Contr., 50 (2005), pp. 841--846..

\bibitem{Bel}
{\sc I. V. Belykh, V. N. Belykh, M. Hasler}, {\em Connection graph
stability method for synchronized coupled chaotic systems},
Physica D, 195 (2004), pp. 159--187.

\bibitem{Sti}
{\sc D. J. Stilwell, E. M. Bollt, D. G. Roberson}, {\em Sufficient
conditions for fast switching synchronization in time varying
network topologies}, SIAM J. Appl. Dyn. Syst., 5:1 (2006), pp.
140--156.


\bibitem{Haj1}
{\sc J. Hajnal}, {\em Weak ergodicity in nonhomogeneous Markov
chains}, Proc. Camb. Phil. Soc., 54 (1958), pp. 233--246.

\bibitem{Haj2}
{\sc J. Hajnal}, {\em The ergodic properties of nonhomogeneous finite
Markov chains}, Proc. Camb. Phil. Soc., 52 (1956), pp. 67--77.

\bibitem{Shen}
{\sc J. Shen}, {\em A geometric approach to ergodic
non-homogeneous Markov chains}, Wavelet Anal. Multi. Meth., LNPAM,
212 (2000), pp. 341--366.


\bibitem{Dau}
{\sc I. Daubechies, J. C. Lagarias}, {\em Sets of matrices all
infinite product of which converge}, Linear Alg. Appl., 161
(1992), pp. 227--263.

\bibitem{Bohl}
{\sc P. Bohl}, {\em \"{U}ber Differentialungleichungen}, J. Reine
Angew. Math., 144 (1913), 284-313.

\bibitem{Barr}
{\sc L. Barreira, Y. B. Pesin}, {\em Lyapunov Exponents and Smooth
Ergodic Theory}, {\rm University Lecture Series}, AMS, Rhode
Island, 2001.

\bibitem{Col}
{\sc F. Colonius, W. Kliemann}, {\em The Lyapunov spectrum of
families of time-varying matrices}, Trans. Americ. Math. Society,
348:11 (1996), pp. 4398--4408.

\bibitem{God}

{\sc C. Godsil, G. Royle}, {\em Algebraic graph theory},
Springer-Verlag, New York, 2001.


\bibitem{Mil}
{\sc J. Milnor}, {\em On the concept of attractor: correction and
remarks}, Commun. Math. Phys., 102 (1985), pp. 517--519.

\bibitem{Vle}
{\sc L. Dieci, Erik S., Van Vleck}, {\em Computation of a few
Lyapunov exponents for continuous and discrete dynamical systems},
Appl. Numer. Math., 17 (1995), pp. 275--291.





\bibitem{Wolf}
{\sc J. Wolfwitz}, {\em Products of indecomposable, aperiodic,
stochastic matrices}, Proc. Amer. Math. Soc., 14:5 (1963), pp.
733--737.

\bibitem{Ose}
{\sc V. I. Oseledec}, {\em A multiplicative ergodic theorem.
Characteristic Ljapunov, exponents of dynamical systems}, Trans.
Moscow Math. Soc., 19 (1968), pp. 197--231.


\bibitem{Atay-PhysD}
{\sc F. M. Atay, T. Biyikoglu, J. Jost}, {\em Network synchronization:
Spectral versus statistical properties.}, Physica D, 224 (2006), pp.
35--41.



\bibitem{Gri}
{\sc G. Gripenberg}, {\em Computing the joint spectral radius},
Linear Alg. Appl., 234 (1996), pp. 43--60.

\bibitem{Zhou1}
{\sc Q. Chen, X. Zhou}, {\em Characterization of joint spectral
radius via trace}, Linear Alg. Appl., 315 (2000), pp. 175--188.

\bibitem{Zhou2}
{\sc X. Zhou}, {\em Estimates for the joint spectral radius},
Appl. Math. Comp., in press, 2006.

\bibitem{Ser}
{\sc D. Serre}, {\em Matrices: Theory and Applications},
Springer-Verlag, New York, 2002.

\bibitem{Arn}
{\sc L. Arnold}, {\em Random Dynamical Systems}, Springer-Verlag,
Heidelberg, 1998.
\bibitem{Mai}
{\sc R. Mainieri}, {\em Zeta function for the Lyapunov exponent of
a product of random matrices}, Phys. Rev. Lett., 68 (1992), pp.
1965--1968.
\bibitem{Bar}
{\sc A.-L. Barabasi, R. Albert}, {\em Emergence of scaling in
random networks}, Science, 286:15 (1999), pp. 509--512.

\end{thebibliography}
\end{document}